\input amstex
\documentstyle{amsppt}
\input epsf.tex

\def\stydate{May 10, 2002}

\chardef\tempcat\catcode`\@ \ifx\undefined\amstexloaded\input
amstex \else\catcode`\@\tempcat\fi \expandafter\ifx\csname
amsppt.sty\endcsname\relax\input amsppt.sty \fi
\let\tempcat\undefined

\immediate\write16{This is LABEL.DEF by A.Degtyarev <\stydate>}
\expandafter\ifx\csname label.def\endcsname\relax\else
  \message{[already loaded]}\endinput\fi
\expandafter\edef\csname label.def\endcsname{%
  \catcode`\noexpand\@\the\catcode`\@\edef\noexpand\styname{LABEL.DEF}%
  \def\expandafter\noexpand\csname label.def\endcsname{\stydate}%
    \toks0{}\toks2{}}
\catcode`\@11
\def\labelmesg@ {LABEL.DEF: }
{\edef\temp{\the\everyjob\W@{\labelmesg@<\stydate>}}
\global\everyjob\expandafter{\temp}}

\def\@car#1#2\@nil{#1}
\def\@cdr#1#2\@nil{#2}
\def\eat@bs{\expandafter\eat@\string}
\def\eat@ii#1#2{}
\def\eat@iii#1#2#3{}
\def\eat@iv#1#2#3#4{}
\def\@DO#1#2\@{\expandafter#1\csname\eat@bs#2\endcsname}
\def\@N#1\@{\csname\eat@bs#1\endcsname}
\def\@Nx{\@DO\noexpand}
\def\@Name#1\@{\if\@undefined#1\@\else\@N#1\@\fi}
\def\@Ndef{\@DO\def}
\def\@Ngdef{\global\@Ndef}
\def\@Nedef{\@DO\edef}
\def\@Nxdef{\global\@Nedef}
\def\@Nlet{\@DO\let}
\def\@undefined#1\@{\@DO\ifx#1\@\relax\true@\else\false@\fi}
\def\@@addto#1#2{{\toks@\expandafter{#1#2}\xdef#1{\the\toks@}}}
\def\@@addparm#1#2{{\toks@\expandafter{#1{##1}#2}%
    \edef#1{\gdef\noexpand#1####1{\the\toks@}}#1}}
\def\make@letter{\edef\t@mpcat{\catcode`\@\the\catcode`\@}\catcode`\@11 }
\def\donext@{\expandafter\egroup\next@}
\def\x@notempty#1{\expandafter\notempty\expandafter{#1}}
\def\lc@def#1#2{\edef#1{#2}%
    \lowercase\expandafter{\expandafter\edef\expandafter#1\expandafter{#1}}}
\newif\iffound@
\def\find@#1\in#2{\found@false
    \DNii@{\ifx\next\@nil\let\next\eat@\else\let\next\nextiv@\fi\next}%
    \edef\nextiii@{#1}\def\nextiv@##1,{%
    \edef\next{##1}\ifx\nextiii@\next\found@true\fi\FN@\nextii@}%
    \expandafter\nextiv@#2,\@nil}
{\let\head\relax\let\specialhead\relax\let\subhead\relax
\let\subsubhead\relax\let\proclaim\relax
\gdef\let@relax{\let\head\relax\let\specialhead\relax\let\subhead\relax
    \let\subsubhead\relax\let\proclaim\relax}}
\newskip\@savsk
\let\@ignorespaces\ignorespaces
\def\@ignorespacesp{\ifhmode
  \ifdim\lastskip>\z@\else\penalty\@M\hskip-1sp%
        \penalty\@M\hskip1sp \fi\fi\@ignorespaces}
\def\ignorespaces{\protect\@ignorespacesp}
\def\@bsphack{\relax\ifmmode\else\@savsk\lastskip
  \ifhmode\edef\@sf{\spacefactor\the\spacefactor}\fi\fi}
\def\@esphack{\relax
  \ifx\penalty@\penalty\else\penalty\@M\fi   
  \ifmmode\else\ifhmode\@sf{}\ifdim\@savsk>\z@\@ignorespacesp\fi\fi\fi}
\let\@frills@\identity@
\let\@txtopt@\identyty@
\newif\if@star
\newif\if@write\@writetrue
\def\@numopt@{\if@star\expandafter\eat@\fi}
\def\checkstar@#1{\DN@{\@writetrue
  \ifx\next*\DN@####1{\@startrue\checkstar@@{#1}}%
      \else\DN@{\@starfalse#1}\fi\next@}\FN@\next@}
\def\checkstar@@#1{\DN@{%
  \ifx\next*\DN@####1{\@writefalse#1}%
      \else\DN@{\@writetrue#1}\fi\next@}\FN@\next@}
\def\checkfrills@#1{\DN@{%
  \ifx\next\nofrills\DN@####1{#1}\def\@frills@####1{####1\nofrills}%
      \else\DN@{#1}\let\@frills@\identity@\fi\next@}\FN@\next@}
\def\checkbrack@#1{\DN@{%
    \ifx\next[\DN@[####1]{\def\@txtopt@########1{####1}#1}%
    \else\DN@{\let\@txtopt@\identity@#1}\fi\next@}\FN@\next@}
\def\check@therstyle#1#2{\bgroup\DN@{#1}\ifx\@txtopt@\identity@\else
        \DNii@##1\@therstyle{}\def\@therstyle{\DN@{#2}\nextii@}%
    \expandafter\expandafter\expandafter\nextii@\@txtopt@\@therstyle.\@therstyle
    \fi\donext@}

\newread\@inputcheck
\def\@input#1{\openin\@inputcheck #1 \ifeof\@inputcheck \W@
  {No file `#1'.}\else\closein\@inputcheck \relax\input #1 \fi}

\def\loadstyle#1{\edef\next{#1}%
    \DN@##1.##2\@nil{\if\notempty{##2}\else\def\next{##1.sty}\fi}%
    \expandafter\next@\next.\@nil\lc@def\next@\next
    \expandafter\ifx\csname\next@\endcsname\relax\input\next\fi}

\let\pagebody@\pagebody
\let\pagetop@\empty
\let\pagebot@\empty
\let\@Xend\empty
\def\pagebody{\pagetop@\pagebody@\pagebot@\@Xend}
\let\@Xclose\empty

\newwrite\@Xmain
\newwrite\@Xsub
\def\W@X{\write\@Xout}
\def\make@Xmain{\global\let\@Xout\@Xmain\global\let\end\endmain@
  \xdef\@Xname{\jobname}\xdef\@inputname{\jobname}}
\begingroup
\catcode`\(\the\catcode`\{\catcode`\{12
\catcode`\)\the\catcode`\}\catcode`\}12
\gdef\W@count#1((\lc@def\@tempa(#1)%
    \def\\##1(\W@X(\global##1\the##1))%
    \edef\@tempa(\W@X(%
        \string\expandafter\gdef\string\csname\space\@tempa\string\endcsname{)%
        \\\pageno\\\cnt@toc\\\cnt@idx\\\cnt@glo\\\footmarkcount@
        \@Xclose\W@X(}))\expandafter)\@tempa)
\endgroup
\def\readaux{\bgroup\checkbrack@\readaux@}
\let\begin\readaux
\def\readaux@{%
    \W@{>>> \labelmesg@ Run this file twice to get x-references right}%
    \global\everypar{}%
    {\def\\{\global\let}%
        \def\/##1##2{\gdef##1{\wrn@command##1##2}}%
        \disablepreambule@cs}%
    \make@Xmain{\make@letter\setboxz@h{\@input{\@txtopt@{\@Xname.aux}}%
            \lc@def\@tempa\jobname\@Name\open@\@tempa\@}}%
  \immediate\openout\@Xout\@Xname.aux%
    \immediate\W@X{\relax}\egroup}
\everypar{\global\everypar{}\readaux}
{\toks@\expandafter{\topmatter}
\global\edef\topmatter{\noexpand\readaux\the\toks@}}
\let\@@end@@\end

\def\@Xclose@{{\def\@Xend{\ifnum\insertpenalties=\z@
        \W@count{close@\@Xname}\closeout\@Xout\fi}%
    \vfill\supereject}}
\def\endmain@{\@Xclose@
    \W@{>>> \labelmesg@ Run this file twice to get x-references right}%
    \@@end@@}
\def\disablepreambule@cs{\\\disablepreambule@cs\relax}

\def\include#1{\bgroup
  \ifx\@Xout\@Xsub\DN@{\errmessage
        {\labelmesg@ Only one level of \string\include\space is supported}}%
    \else\edef\@tempb{#1}\clearpage
      \DN@##1 {\if\notempty{##1}\edef\@tempb{##1}\DN@####1\eat@ {}\fi\next@}%
    \DNii@##1.{\edef\@tempa{##1}\DN@####1\eat@.{}\next@}%
        \expandafter\next@\@tempb\eat@{} \eat@{} %
    \expandafter\nextii@\@tempb.\eat@.%
        \relaxnext@
      \if\x@notempty\@tempa
          \edef\nextii@{\write\@Xmain{%
            \noexpand\string\noexpand\@input{\@tempa.aux}}}\nextii@
        \ifx\undefined\@includelist\found@true\else
                    \find@\@tempa\in\@includelist\fi
            \iffound@\ifx\undefined\@noincllist\found@false\else
                    \find@\@tempb\in\@noincllist\fi\else\found@true\fi
            \iffound@\lc@def\@tempa\@tempa
                \if\@undefined\close@\@tempa\@\else\edef\next@{\@Nx\close@\@tempa\@}\fi
            \else\xdef\@Xname{\@tempa}\xdef\@inputname{\@tempb}%
                \W@count{open@\@Xname}\global\let\@Xout\@Xsub
            \openout\@Xout\@tempa.aux \W@X{\relax}%
            \DN@{\let\end\endinput\@input\@inputname
                    \@Xclose@\make@Xmain}\fi\fi\fi
  \donext@}
\def\includeonly#1{\edef\@includelist{#1}}
\def\noinclude#1{\edef\@noincllist{#1}}

\def\arabicnum#1{\number#1}

\def\Romannum#1{\expandafter\uppercase\expandafter{\romannumeral#1}}
\def\alphnum#1{\ifcase#1\or a\or b\or c\or d\else\@ialph{#1}\fi}
\def\@ialph#1{\ifcase#1\or \or \or \or \or e\or f\or g\or h\or i\or j\or
    k\or l\or m\or n\or o\or p\or q\or r\or s\or t\or u\or v\or w\or x\or y\or
    z\else\fi}
\def\Alphnum#1{\ifcase#1\or A\or B\or C\or D\else\@Ialph{#1}\fi}
\def\@Ialph#1{\ifcase#1\or \or \or \or \or E\or F\or G\or H\or I\or J\or
    K\or L\or M\or N\or O\or P\or Q\or R\or S\or T\or U\or V\or W\or X\or Y\or
    Z\else\fi}

\def\ST@P{step}
\def\ST@LE{style}
\def\N@M{no}
\def\F@NT{font@}
\outer\def\newcounter{\checkbrack@{\expandafter\newcounter@\@txtopt@{{}}}}
{\let\newcount\relax
\gdef\newcounter@#1#2#3{{%
    \toks@@\expandafter{\csname\eat@bs#2\N@M\endcsname}%
    \DN@{\alloc@0\count\countdef\insc@unt}%
    \ifx\@txtopt@\identity@\expandafter\next@\the\toks@@
        \else\if\notempty{#1}\global\@Nlet#2\N@M\@#1\fi\fi
    \@Nxdef\the\eat@bs#2\@{\if\@undefined\the\eat@bs#3\@\else
            \@Nx\the\eat@bs#3\@.\fi\noexpand\arabicnum\the\toks@@}%
  \@Nxdef#2\ST@P\@{}%
  \if\@undefined#3\ST@P\@\else
    \edef\next@{\noexpand\@@addto\@Nx#3\ST@P\@{%
             \global\@Nx#2\N@M\@\z@\@Nx#2\ST@P\@}}\next@\fi
    \expandafter\@@addto\expandafter\@Xclose\expandafter
        {\expandafter\\\the\toks@@}}}}
\outer\def\copycounter#1#2{%
    \@Nxdef#1\N@M\@{\@Nx#2\N@M\@}%
    \@Nxdef#1\ST@P\@{\@Nx#2\ST@P\@}%
    \@Nxdef\the\eat@bs#1\@{\@Nx\the\eat@bs#2\@}}
\outer\def\everystep{\checkstar@\everystep@}
\def\everystep@#1{\if@star\let\next@\gdef\else\let\next@\@@addto\fi
    \@DO\next@#1\ST@P\@}
\def\counterstyle#1{\@Ngdef\the\eat@bs#1\@}
\def\advancecounter#1#2{\@N#1\ST@P\@\global\advance\@N#1\N@M\@#2}
\def\setcounter#1#2{\@N#1\ST@P\@\global\@N#1\N@M\@#2}
\def\counter#1{\refstepcounter#1\printcounter#1}
\def\printcounter#1{\@N\the\eat@bs#1\@}
\def\refcounter#1{\xdef\@lastmark{\printcounter#1}}
\def\stepcounter#1{\advancecounter#1\@ne}
\def\refstepcounter#1{\stepcounter#1\refcounter#1}
\def\savecounter#1{\@Nedef#1@sav\@{\global\@N#1\N@M\@\the\@N#1\N@M\@}}
\def\restorecounter#1{\@Name#1@sav\@}

\def\warning#1#2{\W@{Warning: #1 on input line #2}}
\def\warning@#1{\warning{#1}{\the\inputlineno}}
\def\wrn@@Protect#1#2{\warning@{\string\Protect\string#1\space ignored}}
\def\wrn@@label#1#2{\warning{label `#1' multiply defined}{#2}}
\def\wrn@@ref#1#2{\warning@{label `#1' undefined}}
\def\wrn@@cite#1#2{\warning@{citation `#1' undefined}}
\def\wrn@@command#1#2{\warning@{Preamble command \string#1\space ignored}#2}
\def\wrn@@option#1#2{\warning@{Option \string#1\string#2\space is not supported}}
\def\wrn@@reference#1#2{\W@{Reference `#1' on input line \the\inputlineno}}
\def\wrn@@citation#1#2{\W@{Citation `#1' on input line \the\inputlineno}}
\let\wrn@reference\eat@ii
\let\wrn@citation\eat@ii
\def\nowarning#1{\if\@undefined\wrn@\eat@bs#1\@\wrn@option\nowarning#1\else
        \@Nlet\wrn@\eat@bs#1\@\eat@ii\fi}
\def\printwarning#1{\if\@undefined\wrn@@\eat@bs#1\@\wrn@option\printwarning#1\else
        \@Nlet\wrn@\eat@bs#1\expandafter\@\csname wrn@@\eat@bs#1\endcsname\fi}
\printwarning\Protect \printwarning\label \printwarning\ref
\printwarning\cite \printwarning\command \printwarning\option

{\catcode`\#=12\gdef\@lH{#}}
\def\@@HREF#1{}
\def\@HREF#1#2{\@@HREF{a #1}{\let\@@HREF\eat@#2}\@@HREF{/a}}
\def\@@Hf#1{file:#1} \let\@Hf\@@Hf
\def\@@Hl#1{\if\notempty{#1}\@lH#1\fi} \let\@Hl\@@Hl
\def\@@Hname#1{\@HREF{name="#1"}{}} \let\@Hname\@@Hname
\def\@@Href#1{\@HREF{href="#1"}} \let\@Href\@@Href
\ifx\undefined\pdfoutput
  \csname newcount\endcsname\pdfoutput
\else
  \def\pdflinkattr{attr{/C [0 0.9 0.9]}}
  \let\pdflinkbegin\empty
  \let\pdflinkend\empty
  \def\@pdfHf#1{file {#1}}
  \def\@pdfHl#1{name {#1}}
  \def\@pdfHname#1{\pdfdest name{#1}xyz\relax}
  \def\@pdfHref#1#2{\pdfstartlink \pdflinkattr goto #1\relax
    \pdflinkbegin#2\pdflinkend\pdfendlink}
  \def\@ifpdf#1#2{\ifnum\pdfoutput>\z@\expandafter#1\else\expandafter#2\fi}
  \def\@Hf{\@ifpdf\@pdfHf\@@Hf}
  \def\@Hl{\@ifpdf\@pdfHl\@@Hl}
  \def\@Hname{\@ifpdf\@pdfHname\@@Hname}
  \def\@Href{\@ifpdf\@pdfHref\@@Href}
\fi
\def\@Hr#1#2{\if\notempty{#1}\@Hf{#1}\fi\@Hl{#2}}
\def\@localHref#1{\@Href{\@Hr{}{#1}}}
\def\@countlast#1{\@N#1last\@}
\def\@@countref#1#2{\global\advance#2\@ne
  \@Nxdef#2last\@{\the#2}\@tocHname{#1\@countlast#2}}
\def\@countref#1{\@DO\@@countref#1@HR\@#1}

\def\Href@@#1{\@N\Href@-#1\@}
\def\Href@#1#2{\@N\Href@-#1\@{\@Hl{@#1-#2}}}
\def\Hname@#1{\@N\Hname@-#1\@}
\def\Hlast@#1{\@N\Hlast@-#1\@}
\def\cntref@#1{\global\@DO\advance\cnt@#1\@\@ne
  \@Nxdef\Hlast@-#1\@{\@DO\the\cnt@#1\@}\Hname@{#1}{@#1-\Hlast@{#1}}}
\def\HyperRefs#1{\global\@Nlet\Hlast@-#1\@\empty
  \global\@Nlet\Hname@-#1\@\@Hname
  \global\@Nlet\Href@-#1\@\@Href}
\def\NoHyperRefs#1{\global\@Nlet\Hlast@-#1\@\empty
  \global\@Nlet\Hname@-#1\@\eat@
  \global\@Nlet\Href@-#1\@\eat@}

\HyperRefs{label} {\catcode`\-11
\gdef\@labelref#1{\Hname@-label{r@-#1}}
\gdef\@xHref#1{\Href@-label{\@Hl{r@-#1}}} }
\HyperRefs{toc}
\def\@HR#1{\if\notempty{#1}\string\@HR{\Hlast@{toc}}{#1}\else{}\fi}



\def\bftext{\ifmmode\fam\bffam\else\bf\fi}
\let\@lastmark\empty
\let\@lastlabel\empty
\def\lastmark{\@lastmark}
\let\lastlabel\empty
\let\everylabel\relax
\let\everylabel@\eat@
\let\everyref\relax
\def\newlabel{\bgroup\everylabel\newlabel@}
\def\newlabel@#1#2#3{\if\@undefined\r@-#1\@\else\wrn@label{#1}{#3}\fi
  {\let\protect\noexpand\@Nxdef\r@-#1\@{#2}}\egroup}
\def\w@ref{\bgroup\everyref\w@@ref}
\def\w@@ref#1#2#3#4{%
  \if\@undefined\r@-#1\@{\bftext??}#2{#1}{}\else%
   \@xHref{#1}{\@DO{\expandafter\expandafter#3}\r@-#1\@\@nil}\fi
  #4{#1}{}\egroup}
\def\@@@xref#1{\w@ref{#1}\wrn@ref\@car\wrn@reference}
\def\@xref#1{\rom{\@@@xref{#1}}}
\let\xref\@xref
\def\pageref#1{\w@ref{#1}\wrn@ref\@cdr\wrn@reference}
\def\thepage{\ifnum\pageno<\z@\romannumeral-\pageno\else\number\pageno\fi}
\def\label@{\@bsphack\bgroup\everylabel\label@@}
\def\label@@#1#2{\everylabel@{{#1}{#2}}%
  \@labelref{#2}%
  \let\thepage\relax
  \def\protect{\noexpand\noexpand\noexpand}%
  \edef\@tempa{\edef\noexpand\@lastlabel{#1}%
    \W@X{\string\newlabel{#2}{{\@lastmark}{\thepage}}{\the\inputlineno}}}%
  \expandafter\egroup\@tempa\@esphack}
\def\label#1{\label@{#1}{#1}}
\def\fn@P@{\relaxnext@
    \DN@{\ifx[\next\DN@[####1]{}\else
        \ifx"\next\DN@"####1"{}\else\DN@{}\fi\fi\next@}%
    \FN@\next@}
\def\eat@fn#1{\ifx#1[\expandafter\eat@br\else
  \ifx#1"\expandafter\expandafter\expandafter\eat@qu\fi\fi}
\def\eat@br#1]#2{}
\def\eat@qu#1"#2{}
{\catcode`\~\active\lccode`\~`\@
\lowercase{\global\let\@@P@~\gdef~{\protect\@@P@}}}
\def\Protect@@#1{\def#1{\protect#1}}
\def\disable@special{\let\W@X@\eat@iii\let\label\eat@
    \def\footnotemark{\protect\fn@P@}%
  \let\footnotetext\eat@fn\let\footnote\eat@fn
    \let\refcounter\eat@\let\savecounter\eat@\let\restorecounter\eat@
    \let\advancecounter\eat@ii\let\setcounter\eat@ii
  \let\ifvmode\iffalse\Protect@@\@@@xref\Protect@@\pageref\Protect@@\nofrills
    \Protect@@\\\Protect@@~}
\let\notoctext\identity@
\def\W@X@#1#2#3{\@bsphack{\disable@special\let\notoctext\eat@
    \def\chapter{\protect\chapter@toc}\let\thepage\relax
    \def\protect{\noexpand\noexpand\noexpand}#1%
  \edef\next@{\if\@undefined#2\@\else\write#2{#3}\fi}\expandafter}\next@
    \@esphack}
\newcount\cnt@toc
\def\writeauxline#1#2#3{\W@X@{\cntref@{toc}\let\tocref\@HR}
  \@Xout{\string\@Xline{#1}{#2}{#3}{\thepage}}}
{\let\newwrite\relax
\gdef\@openin#1{\make@letter\@input{\jobname.#1}\t@mpcat}
\gdef\@openout#1{\global\expandafter\newwrite\csname
tf@-#1\endcsname
   \immediate\openout\@N\tf@-#1\@\jobname.#1\relax}}
\def\@@openout#1{\@openout{#1}%
  \@@addto\readaux@{\immediate\closeout\@N\tf@-#1\@}}
\def\auxlinedef#1{\@Ndef\do@-#1\@}
\def\@Xline#1{\if\@undefined\do@-#1\@\expandafter\eat@iii\else
    \@DO\expandafter\do@-#1\@\fi}
\def\beginW@{\bgroup\def\do##1{\catcode`##112 }\dospecials\do\@\do\"
    \catcode`\{\@ne\catcode`\}\tw@\immediate\write\@N}
\def\endW@toc#1#2#3{{\string\tocline{#1}{#2\string\page{#3}}}\egroup}
\def\do@tocline#1{%
    \if\@undefined\tf@-#1\@\expandafter\eat@iii\else
        \beginW@\tf@-#1\@\expandafter\endW@toc\fi
} \auxlinedef{toc}{\do@tocline{toc}}

\let\protect\empty
\def\Protect#1{\if\@undefined#1@P@\@\PROTECT#1\else\wrn@Protect#1\empty\fi}
\def\PROTECT#1{\@Nlet#1@P@\@#1\edef#1{\noexpand\protect\@Nx#1@P@\@}}
\def\pdef#1{\edef#1{\noexpand\protect\@Nx#1@P@\@}\@Ndef#1@P@\@}

\Protect\operatorname \Protect\operatornamewithlimits
\Protect\qopname@ \Protect\qopnamewl@ \Protect\text
\Protect\topsmash \Protect\botsmash \Protect\smash
\Protect\widetilde \Protect\widehat \Protect\thetag
\Protect\therosteritem
\Protect\Cal \Protect\Bbb \Protect\bold \Protect\slanted
\Protect\roman \Protect\italic \Protect\boldkey
\Protect\boldsymbol \Protect\frak \Protect\goth \Protect\dots
\Protect\cong \Protect\lbrace \let\{\lbrace \Protect\rbrace
\let\}\rbrace
\let\root@P@@\root \def\root@P@#1{\root@P@@#1\of}
\def\root#1\of{\protect\root@P@{#1}}

\def\frills{\ignorespaces\@txtopt@}
\def\frillsnotempty#1{\x@notempty{\@txtopt@{#1}}}
\def\numberline{\@numopt@}
\newif\if@theorem
\let\@therstyle\eat@
\def\@headtext@#1#2{{\disable@special\let\protect\noexpand
    \def\chapter{\protect\chapter@rh}%
    \edef\next@{\noexpand\@frills@\noexpand#1{#2}}\expandafter}\next@}
\let\AmSrighthead@\rightheadtext
\def\rightheadtext{\checkfrills@{\@headtext@\AmSrighthead@}}
\let\AmSlefthead@\leftheadtext
\def\leftheadtext{\checkfrills@{\@headtext@\AmSlefthead@}}
\def\@head@@#1#2#3#4#5{\@Name\pre\eat@bs#1\@\if@theorem\else
    \@frills@{\csname\expandafter\eat@iv\string#4\endcsname}\relax
        \ifx\protect\empty\@N#1\F@NT\@\fi\fi
    \@N#1\ST@LE\@{\counter#3}{#5}%
  \if@write\writeauxline{toc}{\eat@bs#1}{#2{\counter#3}\@HR{#5}}\fi
    \if@theorem\else\expandafter#4\fi
    \ifx#4\endhead\ifx\@txtopt@\identity@\else
        \headmark{\@N#1\ST@LE\@{\counter#3}{\frills\empty}}\fi\fi
    \@Name\post\eat@bs#1\@\ignorespaces}
\ifx\undefined\endhead\Invalid@\endhead\fi
\def\@head@#1{\checkstar@{\checkfrills@{\checkbrack@{\@head@@#1}}}}
\def\@thm@@#1#2#3{\@Name\pre\eat@bs#1\@
    \@frills@{\csname\expandafter\eat@iv\string#3\endcsname}
    {\@theoremtrue\check@therstyle{\@N#1\ST@LE\@}\frills
            {\counter#2}\@theoremfalse}%
    \@DO\envir@stack\end\eat@bs#1\@
    \@N#1\F@NT\@\@Name\post\eat@bs#1\@\ignorespaces}
\def\@thm@#1{\checkstar@{\checkfrills@{\checkbrack@{\@thm@@#1}}}}
\def\@capt@@#1#2#3#4#5\endcaption{\bgroup
    \edef\@tempb{\global\footmarkcount@\the\footmarkcount@
    \global\@N#2\N@M\@\the\@N#2\N@M\@}%
    \def\shortcaption##1{\global\def\sh@rtt@xt####1{##1}}\let\sh@rtt@xt\identity@
    \DN@{#4{\@tempb\@N#1\ST@LE\@{\counter#2}}}%
    \if\notempty{#5}\DNii@{\next@\@N#1\F@NT\@}\else\let\nextii@\next@\fi
    \nextii@#5\endcaption
  \if@write\writeauxline{#3}{\eat@bs#1}{{} \@HR{\@N#1\ST@LE\@{\counter#2}%
    \if\notempty{#5}.\enspace\fi\sh@rtt@xt{#5}}}\fi
  \global\let\sh@rtt@xt\undefined\egroup}
\def\@capt@#1{\checkstar@{\checkfrills@{\checkbrack@{\@capt@@#1}}}}
\let\captiontextfont@\empty

\ifx\undefined\subsubheadfont@\def\subsubheadfont@{\it}\fi
\ifx\undefined\proclaimfont\def\proclaimfont{\sl}\fi
\ifx\undefined\proclaimfont@\let\proclaimfont@\proclaimfont\fi
\def\proclaimfont{\proclaimfont@}
\ifx\undefined\definitionfont@\def\AmSdeffont@{\rm}
    \else\let\AmSdeffont@\definitionfont@\fi
\ifx\undefined\remarkfont@\def\remarkfont@{\rm}\fi

\def\newfont@def#1#2{\if\@undefined#1\F@NT\@
    \@Nxdef#1\F@NT\@{\@Nx.\expandafter\eat@iv\string#2\F@NT\@}\fi}
\def\newhead@#1#2#3#4{{%
    \gdef#1{\@therstyle\@therstyle\@head@{#1#2#3#4}}\newfont@def#1#4%
    \if\@undefined#1\ST@LE\@\@Ngdef#1\ST@LE\@{\headstyle}\fi
    \if\@undefined#2\@\gdef#2{\headtocstyle}\fi
  \@@addto\moretocdefs@{\\#1#1#4}}}
\outer\def\newhead#1{\checkbrack@{\expandafter\newhead@\expandafter
    #1\@txtopt@\headtocstyle}}
\outer\def\newtheorem#1#2#3#4{{%
    \gdef#2{\@thm@{#2#3#4}}\newfont@def#2#4%
    \@Nxdef\end\eat@bs#2\@{\noexpand\revert@envir
        \@Nx\end\eat@bs#2\@\noexpand#4}%
  \if\@undefined#2\ST@LE\@\@Ngdef#2\ST@LE\@{\proclaimstyle{#1}}\fi}}%
\outer\def\newcaption#1#2#3#4#5{{\let#2\relax
  \edef\@tempa{\gdef#2####1\@Nx\end\eat@bs#2\@}%
    \@tempa{\@capt@{#2#3{#4}#5}##1\endcaption}\newfont@def#2\endcaptiontext%
  \if\@undefined#2\ST@LE\@\@Ngdef#2\ST@LE\@{\captionstyle{#1}}\fi
  \@@addto\moretocdefs@{\\#2#2\endcaption}\newtoc{#4}}}
{
\outer\gdef\newtoc#1{{%
    \@DO\ifx\do@-#1\@\relax
    \global\auxlinedef{#1}{\do@tocline{#1}}{}%
    \@@addto\tocsections@{\make@toc{#1}{}}\fi}}}

\toks@\expandafter{\itembox@}
\toks@@{\bgroup\let\therosteritem\identity@\let\rm\empty
  \let\@Href\eat@\let\@Hname\eat@
  \edef\next@{\edef\noexpand\@lastmark{\therosteritem@}}\donext@}
\edef\itembox@{\the\toks@@\the\toks@}
\def\firstitem@false{\let\iffirstitem@\iffalse
    \global\let\lastlabel\@lastlabel}

\let\rosteritemrefform\therosteritem
\let\rosteritemrefseparator\empty
\def\rosteritemref#1{\hbox{\rosteritemrefform{\@@@xref{#1}}}}
\def\local#1{\label@\@lastlabel{\lastlabel-i#1}}

\def\xRef@P@{\gdef\lastlabel}
\def\xRef#1{\@xref{#1}\protect\xRef@P@{#1}}

\def\iref@P@{\gdef\lastref}
\def\itemref#1#2{\rosteritemref{#1-i#2}\protect\iref@P@{#1}}
\def\iref#1{\@xref{#1}\rosteritemrefseparator\itemref{#1}}

\def\eqref#1{\thetag{\@@@xref{#1}}}
\def\tagform@#1{\ifmmode\hbox{\rm\else\rom{\fi
        (\ignorespaces#1\unskip)\iftrue}\else}\fi}

\let\AmSfnote@\makefootnote@
\def\makefootnote@#1{\bgroup\let\footmarkform@\identity@
  \edef\next@{\edef\noexpand\@lastmark{#1}}\donext@\AmSfnote@{#1}}

\def\clearpage{\ifnum\insertpenalties>0\line{}\fi\vfill\supereject}

\def\proof{\checkfrills@{\checkbrack@{%
    \check@therstyle{\@frills@{\demo}{\frills{Proof}}{}}
        {\frills{}\envir@stack\endremark\envir@stack\enddemo}%
  \envir@stack\endproof\ignorespaces}}}
\def\everyendproof{\qed}
\def\endproof{\nofrillscheck{\frills@{\everyendproof}\revert@envir\endproof\enddemo}}

\let\AmSref\ref
\let\AmSrefstyle\refstyle
\let\plaincite\cite
\def\citei@#1,{\citeii@#1\eat@,}
\def\citeii@#1\eat@{\w@ref{#1}\wrn@cite\@car\wrn@citation}
\def\mcite@#1;{\plaincite{\citei@#1\eat@,\unskip}\mcite@i}
\def\mcite@i#1;{\DN@{#1}\ifx\next@\endmcite@
  \else, \plaincite{\citei@#1\eat@,\unskip}\expandafter\mcite@i\fi}
\def\endmcite@{\endmcite@}
\def\cite#1{\mcite@#1;\endmcite@;}
\PROTECT\cite
\def\refstyle#1{\AmSrefstyle{#1}\uppercase{%
    \ifx#1A\relax \def\@ref@##1{\AmSref\xdef\@lastmark{##1}\key##1}%
    \else\ifx#1C\relax \def\@ref@##1{\AmSref\no\counter\refno}%
        \else\def\@ref@{\AmSref}\fi\fi}}
\refstyle A
\newcounter\refno\null
\newif\ifRefs
\gdef\Refs{\checkstar@{\checkbrack@{\csname AmSRefs\endcsname
  \nofrills{\frills{References}%
  \if@write\writeauxline{toc}{vartocline}{\@HR{\frills{References}}}\fi}%
  \def\ref{\@ref@}\Refstrue\ignorespaces}}}
\let\ref\xref

\newif\iftoc
\pdef\tocbreak{\iftoc\hfil\break\fi}
\def\tocsections@{\make@toc{toc}{}}
\let\moretocdefs@\empty
\def\newtocline@#1#2#3{%
  \edef#1{\def\@Nx#2line\@####1{\@Nx.\expandafter\eat@iv
        \string#3\@####1\noexpand#3}}%
  \@Nedef\no\eat@bs#1\@{\let\@Nx#2line\@\noexpand\eat@}%
    \@N\no\eat@bs#1\@}
\def\MakeToc#1{\@@openout{#1}}
\def\newtocline#1#2#3{\Err@{\Invalid@@\string\newtocline}}
\def\make@toc#1#2{\penaltyandskip@{-200}\aboveheadskip
    \if\notempty{#2}
        \centerline{\headfont@\ignorespaces#2\unskip}\nobreak
    \vskip\belowheadskip \fi
    \@openin{#1}\relax
    \vskip\z@}
\def\contents{\readaux\checkfrills@{\checkbrack@{\@contents@}}}
\def\@contents@{\toc@{\frills{Contents}}\envir@stack\endcontents%
    \def\nopagenumbers{\let\page\eat@}\let\newtocline\newtocline@\toctrue
  \def\@HR{\Href@{toc}}%
  \def\tocline##1{\csname##1line\endcsname}
  \edef\caption##1\endcaption{\expandafter\noexpand
    \csname head\endcsname##1\noexpand\endhead}%
    \ifmonograph@\def\vartoclineline{\Chapterline}%
        \else\def\vartoclineline##1{\sectionline{{} ##1}}\fi
  \let\\\newtocline@\moretocdefs@
    \ifx\@frills@\identity@\def\\##1##2##3{##1}\moretocdefs@
        \else\let\tocsections@\relax\fi
    \def\\{\unskip\space\ignorespaces}\let\maketoc\make@toc}
\def\endcontents{\tocsections@\vskip-\lastskip\revert@envir\endcontents
    \endtoc}

\if\@undefined\selectf@nt\@\let\selectf@nt\identity@\fi
\def\Err@math#1{\Err@{Use \string#1\space only in text}}
\def\textonlyfont@#1#2{%
    \def#1{\RIfM@\Err@math#1\else\edef\f@ntsh@pe{\string#1}\selectf@nt#2\fi}%
    \PROTECT#1}
\tenpoint

\def\newshapeswitch#1#2{\gdef#1{\selectsh@pe#1#2}\PROTECT#1}
\def\shapeswitch#1#2#3{\@Ngdef#1\string#2\@{#3}}
\shapeswitch\rm\bf\bf \shapeswitch\rm\tt\tt
\shapeswitch\rm\smc\smc
\newshapeswitch\em\it
\shapeswitch\em\it\rm \shapeswitch\em\sl\rm
\def\selectsh@pe#1#2{\relax\if\@undefined#1\f@ntsh@pe\@#2\else
    \@N#1\f@ntsh@pe\@\fi}

\def\@itcorr@{\leavevmode
    \edef\prevskip@{\ifdim\lastskip=\z@ \else\hskip\the\lastskip\relax\fi}\unskip
    \edef\prevpenalty@{\ifnum\lastpenalty=\z@ \else
        \penalty\the\lastpenalty\relax\fi}\unpenalty
    \/\prevpenalty@\prevskip@}
\def\rom@P@#1{\@itcorr@{\selectsh@pe\rm\rm#1}}
\def\rom{\protect\rom@P@}
\def\Rom@P@#1{\@itcorr@{\rm#1}}
\def\Rom{\protect\Rom@P@}
{\catcode`\-11 \HyperRefs{idx} \HyperRefs{glo}
\newcount\cnt@idx \global\cnt@idx=10000
\newcount\cnt@glo \global\cnt@glo=10000
\gdef\writeindex#1{\W@X@{\cntref@{idx}}\tf@-idx
 {\string\indexentry{#1}{\Hlast@{idx}}{\thepage}}}
\gdef\writeglossary#1{\W@X@{\cntref@{glo}}\tf@-glo
 {\string\glossaryentry{#1}{\Hlast@{glo}}{\thepage}}}
}
\def\emph#1{\@itcorr@\bgroup\em\ignorespaces#1\unskip\egroup
  \DN@{\DN@{}\ifx\next.\else\ifx\next,\else\DN@{\/}\fi\fi\next@}\FN@\next@}
\def\makequoteactive{\catcode`\"\active}
{\makequoteactive\gdef"{\FN@\quote@}
\gdef\quote@{\ifx"\next\DN@"##1""{\quoteii{##1}}\else\DN@##1"{\quotei{##1}}\fi\next@}}
\let\quotei\eat@
\let\quoteii\eat@
\def\MakeIndex{\@openout{idx}}
\def\MakeGlossary{\@openout{glo}}

\def\endofpar#1{\ifmmode\ifinner\endofpar@{#1}\else\eqno{#1}\fi
    \else\leavevmode\endofpar@{#1}\fi}
\def\endofpar@#1{\unskip\penalty\z@\null\hfil\hbox{#1}\hfilneg\penalty\@M}

\newdimen\normalparindent\normalparindent\parindent
\def\firstparindent#1{\everypar\expandafter{\the\everypar
  \global\parindent\normalparindent\global\everypar{}}\parindent#1\relax}

\@@addto\disablepreambule@cs{%
    \\\readaux\relax
    \\\begin\relax
    \\\readaux@\relax
    \\\@openout\eat@
    \\\@@openout\eat@
    \/\Monograph\empty
    \/\MakeIndex\empty
    \/\MakeGlossary\empty
    \/\MakeToc\eat@
    \/\HyperRefs\eat@
    \/\NoHyperRefs\eat@
}

\csname label.def\endcsname


\def\punct#1#2{\if\notempty{#2}#1\fi}
\def\sppunct{\punct{.\enspace}}
\def\varpunct#1#2{\if\frillsnotempty{#2}#1\fi}

\def\headstyle#1#2{\numberline{#1\sppunct{#2}}\ignorespaces#2\unskip}
\def\headtocstyle#1#2{\numberline{#1\punct.{#2}}\space #2}

\def\specialtocstyle#1#2{#2}
\newcounter\section\null
\newcounter\subsection\section
\newcounter\subsubsection\subsection
\newhead\specialsection[\specialtocstyle]\null\endspecialhead
\newhead\section\section\endhead
\newhead\subsection\subsection\endsubhead
\newhead\subsubsection\subsubsection\endsubsubhead
\def\firstappendix{\global\sectionno0 %
  \counterstyle\section{\Alphnum\sectionno}%
    \global\let\firstappendix\empty}

\def\appendixtocstyle#1#2{\space\numberline{Appendix #1\sppunct{#2}}#2}
\newhead\appendix[\appendixtocstyle]\section\endhead

\let\endAmSdef\enddefinition
\def\proclaimstyle#1#2{\numberline{#2\varpunct{.\enspace}{#1}}\frills{#1}}
\copycounter\thm\subsubsection
\theorem\thm\endproclaim
\proposition\thm\endproclaim
\lemma\thm\endproclaim
\corollary\thm\endproclaim
\definition\thm\endAmSdef
\example\thm\endAmSdef

\def\captionstyle#1#2{\frills{#1}\numberline{\varpunct{ }{#1}#2}}
\newcounter\figure\null
\newcounter\table\null
\newcaption{Figure}\figure\figure{lof}\botcaption
\newcaption{Table}\table\table{lot}\topcaption

\copycounter\equation\subsubsection


\font\sani=cmssi10
\font\fib=cmfib8

\font\tenscr=rsfs10 \font\sevenscr=rsfs7 \font\fivescr=rsfs5
\skewchar\tenscr='177 \skewchar\sevenscr='177
\skewchar\fivescr='177
\newfam\scrfam \textfont\scrfam=\tenscr \scriptfont\scrfam=\sevenscr
\scriptscriptfont\scrfam=\fivescr
\define\scr#1{{\fam\scrfam#1}}

\pagewidth{13.3cm}
\def\BB{\scr B}
\def\CC{\scr C}
\def\DD{\scr D}
\def\conj{\operatorname{conj}}
\def\C{{\Bbb C}}
\def\R{{\Bbb R}}
\def\Z{{\Bbb Z}}

\def\Rp#1{\R\roman P^{#1}}
\def\Cp#1{\C \roman P^{#1}}
\def\Im{\mathop{\roman{Im}}\nolimits}
\def\Re{\mathop{\roman{Re}}\nolimits}

\def\ind{\mathop{\roman{ind}}\nolimits}
\def\Int{\mathop{\roman{Int}}\nolimits}
\def\Cl{\mathop{\roman{Cl}}\nolimits}
\def\la{\langle}
\def\ra{\rangle}

\def\Aut{\mathop{\roman{Aut}}\nolimits}

\def\oo{\varnothing}

\def\per{\operatorname{per}}
\def\rank{\operatorname{rank}}

\def\D{\Delta}

\def\til{\widetilde}

\def\e{\varepsilon}
\def\a{\alpha}
\def\D{\Delta}

\def\Cl{\mathop{\roman{Cl}}\nolimits}
\def\PP{P}
\def\QQ{Q}
\def\QQQ{\Cal Q}
\def\VV{V}
\def\LL{L}
\def\ZZ{Z}

\def\Cr{C_{\text{root}}}
\def\mm{\text{\sani m}}
\def\sm{\smallsetminus}
\def\pt{\roman{pt}}

\def\M{\Bbb M}
\def\LLL{\Bbb L}

\let\rk=\remark
\let\endrk=\endremark
\let\ge\geqslant
\let\le\leqslant
\let\+\sqcup
\let\dsum\+


\NoBlackBoxes

 \topmatter
\title
Topology of real cubic fourfolds
\endtitle
\author S.~Finashin, V.~Kharlamov
\endauthor
\address Middle East Technical University,
Department of Mathematics\endgraf Ankara 06531 Turkey
\endaddress
\address
Universit\'{e} de Strasbourg et IRMA (CNRS)\endgraf 7 rue Ren\'{e}
Descartes 67084 Strasbourg Cedex, France
\endaddress
 \subjclass 14P25, 14J10, 14N25, 14J35,
14J70
\endsubjclass
\abstract A solution of the problem of topological classification
of real cubic fourfolds is given. It is proven that the real locus
of a real non-singular cubic fourfold is diffeomorphic either to a
connected sum $\Rp4\#i(S^2\times S^2)\#j(S^1\times S^3)$, or to a
disjoint union $\Rp4\+ S^4$.
\endabstract
\endtopmatter

\document

\rightline{\vbox{\hsize 75mm \noindent\eightit\baselineskip10pt La
voie la plus courte et la meilleure entre deux v\'erit\'es de
domaine r\'eel passe souvent par le domaine imaginaire.
 \vskip3mm\noindent\eightrm
 J.S. Hadamard, Essai sur la psychologie de l'invention dans le domaine math\'ematique, Edition Gaithier-Villars, 1975, p. 114}}
 \vskip10mm

\section {Introduction}

\subsection{The subject}
Studying the topology of cubic hypersurfaces is a classical task
in real algebraic geometry. It has a long history going back to
Newton, who undertook a systematic study of real plane cubics
(first published in 1704 as an appendix to his book {\it
Opticks}), and to Schl\"{a}fli, Cayley, and Klein, who in a series
of treatises (dating from 1858 to 1873 and certainly strongly
motivated by Cayley-Salmon discovery of 27 straight lines on
nonsingular cubic surfaces) classified the shapes of real cubic
surfaces. Cubic surfaces (as well as cubic curves and higher
dimensional cubic hypersurfaces) can be classified in different
manners. As far as we know, it was Klein who first clearly
addressed the problem of deformation and topological
classifications. He solved it for nonsingular cubic surfaces and
specially emphasized that for nonsingular cubic surfaces the
deformation classification coincides with the topological one: two
real nonsingular cubic surfaces in $P^3$ are deformation
equivalent if and only if their real point sets are homeomorphic.
Namely, there are $5$ deformation classes of nonsingular cubic
surfaces respectively to $5$ topological types: $\Rp2 \+S^2$,
$\Rp2$, $\#_3\Rp2$, $\#_5\Rp2$, and $\#_7\Rp2$ ($\#$ stands for
the connected sum and $\+$, for the disjoint sum). Recall, that
the situation with the real plane cubics is similar: there are $2$
deformation classes of nonsingular cubic curves and they are
distinguished by their topological types, $S^1 \+S^1$ and $S^1$
(or one may like to write more instructively, $\Rp1\+S^1$ and
$\Rp1$).

Only recently a deformation classification of real cubic
threefolds was completed by V.~Krasnov~\cite{Kr1} who proved that
there are $9$ classes of real nonsingular cubic threefolds
$X\subset P^4$. It turned out that these $9$ classes of cubics $X$
are distinguished by the Betti numbers of the real point set
$X(\R)$ plus vanishing or non-vanishing of the homology class
realized by $X(\R)$ in the middle homology group $H_3(X;\Z/2)$ of
the complex point set (we denote the complex point set by the same
letter as the variety itself). In a subsequent paper~\cite{Kr2}
Krasnov determined the topological type of $X(\R)$ for $8$ of the
9 deformation classes: for 1 class it is $\Rp3 \+S^3$ and for $7$
others it is $\Rp3\#_k(S^1\times S^2)$, $k=0,\dots,6$.

In \cite{FK1, FK2} we have undertaken a systematic study of
deformations classes of real cubic fourfolds. At the first step,
in \cite{FK1} we obtained a {\it coarse deformation
classification} (that is a classification up to deformations
combined with the projective equivalence) of real cubic fourfolds.
There we did not only enumerate the deformation classes (their
number is 75), but also described them in terms of simple
homological invariants, and, in addition, provided the adjacency
graph (which we call {\it the $\roman K4$-graph})
 of the deformation classes. This graph essentially
coincides with the adjacency graph for non-polarized real
K3-surfaces ({\it the $\roman K3$-graph}), which was a somewhat
unexpected outcome of our research, even though cubic fourfolds
are well-known to be related to the K3 surfaces in many ways.

The difference between the coarse and ordinary deformation
classifications is encoded in the chirality phenomena: we say that
a real nonsingular cubic $X\subset P^5$ and its coarse deformation
class are {\it chiral} if $X$ and its image under a mirror
reflection, $X'$, belong to different connected components of the
space of real nonsingular cubics. Thus, each chiral coarse
deformation class gives two ordinary deformation classes, and each
achiral coarse deformation class gives only one ordinary
deformation class. In \cite{FK2} we analyzed the chirality of the
deformation classes of real nonsingular cubic fourfolds. We
reduced the chirality problem to a specific problem of the
arithmetics of lattices and used this reduction to show that
certain real cubic fourfolds are chiral, while certain other real
cubic fourfolds are achiral. The crucial role in the reduction to
the arithmetics of lattices is played by the period map and the
corresponding surjectivity statement (which was established in the
complex setting by R.~Laza \cite{Laza} and E.~Looijenga
\cite{Looijenga} shortly before our work).

\subsection{The main result}
Despite such a detailed understanding of the deformation classes,
the topology of the real locus for real nonsingular cubic
fourfolds being too much beyond the control under indicated above
approaches remained far from being understood. The topological
classification required additional tools.

In the present paper we introduce two such ingredients. One of
them is based on detecting cuspidal cubics on the boundary of
deformation components and analyzing the surgeries provided by
perturbations of a cuspidal cubic. Another essential ingredient is
a technique of {\it ramified connected sums}. Combining these
tools with the results and certain methods from \cite{FK1, FK2},
we solve the problem of topological classification of real cubic
fourfolds. Our main result is the following (more detailed
statements are given in Section \ref{main-thm}).

\theorem\label{main-theorem} The real locus of a real non-singular
cubic fourfold is diffeomorphic either to a connected sum
$\Rp4\#i(S^2\times S^2)\#j(S^1\times S^3)$ {\rm(}the connected
case{\rm)}, or to a disjoint union of $\Rp4$ with $S^4$
{\rm(}disconnected case{\rm)}.
\endtheorem

In addition, we complete the topological description of real
nonsingular cubic threefolds, $X\subset P^4$, by applying our
technique to the deformation class which was not analyzed by
Krasnov. In this case the real locus, $X(\R)$, turns out to be a
Seifert manifold described in Section~\ref{remarks} (Theorem
\ref{spiral}).

\subsection{Structure of the paper}
The paper is organized as follows.

In Section \ref{chambers} we review some basic facts about real
cubic fourfolds. In particular, in Subsection \ref{K4} we describe
the K4-graph, whose vertices $C^{i,j}$ and $C^{i,j}_I$ represent
the deformation classes of real cubic fourfolds and whose edges
describe adjacency of these classes. The first step towards
proving Theorem \ref{main-theorem} is made in Subsection
\ref{immediate}, where using elementary arguments, we deduce that
vertices $C^{0,0}$ and $C^{1,0}_I$ represent cubic fourfolds $X$
whose real locus $X(\R)$ is $\Rp4$ and $\Rp4\+S^4$ respectively.
Subsections 2.6---2.7 contain basic definitions and facts related
to the period map for cubic fourfolds.
 This techniques (which played a crucial role in \cite{FK2})
is involved here only in one place: in the proof of Lemma
\ref{cusp-lemma} (this Lemma provides an arithmetical criterion of
existence of cuspidal cubics on the boundary of a given
deformation component).

In Section \ref{AiT} we analyze the Morse modifications
experienced by $X(\R)$ as we move down the K4-graph. This
technique allows us to use some kind of inductive scheme for
proving Theorem \ref{main-theorem} with $C^{0,0}$ as the basis of
induction. However it turns out that in two cases,
 $C^{2,1}_I$ and $C^{10,1}$, such
Morse theory inductive arguments do not work, and we need to apply
a more involved technique. Such a technique of "ramified connected
sums" is developed in Section \ref{ramified}.

In Section \ref{main-thm} we apply the results and the technique
of the previous sections to prove Theorem \ref{main-theorem}. In
Section \ref{remarks} we make some concluding remarks. In
particular, we apply the technique of ramified connected sums to
complete the topological classification of real cubic threefolds
(Theorem \ref{spiral}). We discuss also certain perspectives in
studying higher dimensional cubics.

\subsection{Conventions}
When in a subsection all or most of the material is a recollection
of already known results, the appropriate references are given in
brackets after the subsection title (some results may be not
explicitly mentioned in those cited papers, in such a case we
provide their proof). As usual, symbol `$\square$' after a
statement means that no proof will follow: either the proof is
straightforward, or it has already been explained, or a reference
is given (in the statement, or at the beginning of the
subsection).

Notation $X(\R)$ refers always to the real locus of a variety $X$
defined over $\R$. However we write $\Rp{n}$ instead of $P^n(\R)$
in the topological context. In Section 4, we also diverge from
algebro-geometric notation style by using single letters ($P$,
$Q$, $V$, $L$) for certain real loci as soon as they are involved
in topological constructions.

\subsection{Acknowledgements}
The principal ideas of this paper were developed by the authors
during their visit to Centre Interfacultaire Bernoulli in EPFL
(Lausanne). The text was finalized during the first author's visit
to Universit\'e de Strasbourg. The last touch was made during
RIP-stay in Mathematiches Forschungsinstitut Oberwolfach. We thank
these institutions for hospitality.

\section{Deformation chambers and the period map}
\label{chambers}

\subsection{Discriminant hypersurface and  deformation
chambers (\cite{FK1})} The complex projective cubic fourfolds form
the complex projective space $P_{4,3}=P(Sym^3 (\C^6)^*)$ of
dimension $\binom{5+3}3-1=55$ and the singular ones form in it
the, so called, \emph{discriminant hypersurface} $\Delta\subset
P_{4,3}$. The real cubic fourfolds, that is cubic fourfolds
defined by real polynomials, form the real projective space
$P_{4,3}(\R)$, so that the real singular cubic fourfolds form {\it
the real discriminant hypersurface} $\Delta(\R)\subset
P_{4,3}(\R)$.

We study the non-singular real cubics, so the space of our
interest is nothing but the complement $\CC=P_{4,3}(\R)
\sm\Delta(\R)$ of the discriminant hypersurface. We say that
non-singular real cubic fourfolds $X_1$ and $X_2$ are {\it
deformation equivalent} if they belong to the same component of
$\CC$. The connected components of $\CC$ will be called {\it the
deformation components}, or {\it the deformation chambers}.

\remark{{\eightit Remark}}{\eightrm In our definition we consider
a real cubic as a singular one, even if it has only imaginary
singular points. Thus, it may be worth mentioning that the real
cubics with real singular points form a real semi-algebraic subset
$\Delta'\subset\Delta$, such that the difference $\Delta \sm
\Delta'$ has codimension $2$ in $P_{4,3}(\R)$. Because of this
phenomenon the imbedding of $P_{4,3}(\R)\sm\D$ into
$P_{4,3}(\R)\sm\D'$ induces a bijection at the level of connected
components and the two corresponding deformation classifications
coincide.}
\endremark\vskip0.05in

The smooth part $\Delta_1$ of $\Delta$ is formed by cubics with a
node ($A_1$-singularity) and no other singular points. The
connected components of $\Delta_1(\R)$ are called {\it facets}.
Two facets can be adjacent through a stratum formed by cubics with
a cusp ($A_2$-singularity), in such a case the two facets are
called {\it cuspidally adjacent}. The union of the facets and the
strata formed by cubics with a cusp is a topological
codimension-one submanifold of $P_{4,3}(\R)$. We denote this union
by $\D_2(\R)$ and call its connected components the {\it walls}.
 If two deformation components $C$ and $C'$ have a common wall,  we say that
they are {\it adjacent}.

A cubic $X\in\CC$ is called {\it achiral} if it is deformation
equivalent to its image $\rho(X)$ under the mirror reflection
$\rho\:P^5(\R)\to P^5(\R)$ against some hyperplane in $P^5(\R)$;
otherwise the cubic $X$ is called {\it chiral}.

The relation of {\it coarse deformation equivalence} is a
combination of deformation and projective equivalences. Namely,
cubics $X_1$ and $X_2$ are {\it coarse deformation equivalent} if
$X_1$ is deformation equivalent to the image of $X_2$ under a
projective transformation (that is either to $X_2$ or to
$\rho(X_2)$).

In accord with this, by {\it the coarse deformation component} we
mean a connected component of the quotient-space
$\til{\CC}=\CC/PGL(6,\R)$ with respect to the action of the
projective group $PGL(6,\R)$ in $P_{4,3}(\R)$ induced by the
action in $P^5(\R)$.

It is obvious that one coarse deformation component of chiral
cubics corresponds to two ordinary deformation components, and
that for achiral cubics the deformation and the coarse deformation
components are in one-to-one correspondence.

A pair of coarse deformation components are called {\it adjacent}
if there is a pair of adjacent deformation components representing
these coarse components.

\subsection{K4- and K3-graphs(\cite{FK1})}\label{K4} Presenting the coarse deformation
components of cubic fourfolds as vertices and their adjacency as
edges, we obtain a graph. We call it the K4{\it-graph} and denote
it by $\Gamma_{K4}$. This name is essentially to reflect its
remarkable similarity to the better known K3{\it-graph} which
describes similarly the set of deformation components of real
non-polarized K3-surfaces and their adjacency.

\midinsert \topcaption{Figure 1. The K4-graph $\Gamma_{K4}$}
\endcaption\epsfbox{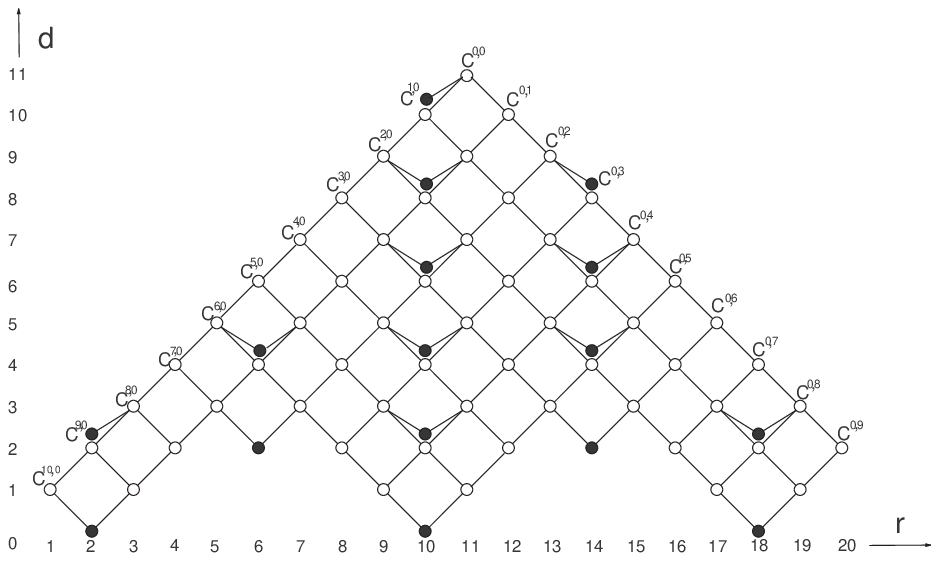}
\eightpoint \noindent
\endinsert

\midinsert \topcaption{Figure 2. The K3-graph $\Gamma_{K3}$}
\endcaption\epsfbox{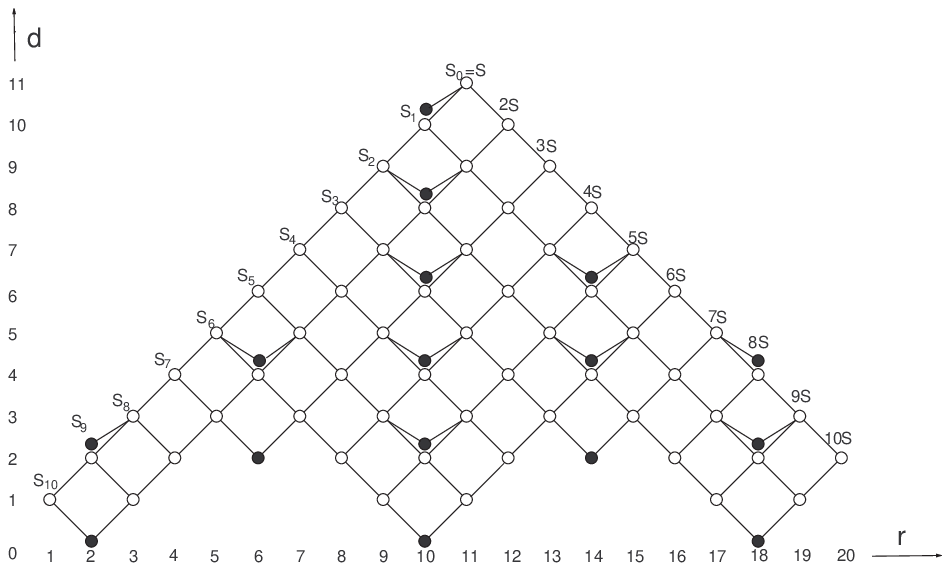}
\eightpoint \noindent
\endinsert

\theorem The graphs $\Gamma_{K4}$ and $\Gamma_{K3}$ are the graphs
shown respectively on Figures 1 and 2. \qed
\endtheorem

The coordinates  $r$ and $d$ on the Figures characterize the
action of the complex conjugation in the corresponding lattices.
Namely, in the case of a real cubic fourfold $X$, we consider the
involution $c_X\:\M\to\M$ induced in $\M=H^4(X)$ by the complex
conjugation, and denote by $d$ the rank of the 2-periodic {\it
discriminant group} $\M/(\M_++\M_-)=\M_\pm^*/\M_\pm$, where
$\M_\pm=\{x\in\M\,|\,c_X(x)=\pm x\}$. By well known Smith theory
arguments, $d=\frac12\big(b_*(X)-b_*(X(\R))\big)=\frac12\big(27
-b_*(X(\R))\big)$ where $b_*=\sum_{i\ge0}b_i$ is the total Betti
number with $\Z/2$-coefficients.
 In the case of a real K3-surface $Y$, there is a similar involution
$c_Y\:\LLL\to\LLL$ in $\LLL=H^2(Y)$ and $d$ denotes the rank of
the discriminant group $\LLL/(\LLL_++\LLL_-)=\LLL_\pm^*/\LLL_\pm$,
where $\LLL_\pm=\{x\in\LLL\,|\,c_Y(x)=\pm x\}$. If $Y(\R)\ne\oo$,
then we have similarly
$d=\frac12\big(b_*(Y)-b_*(Y(\R))\big)=\frac12\big(22-b_*(Y(\R))\big)$,
and in the exceptional case $Y(\R)=\oo$, we have $d=10$.

 The coordinate $r$ is the rank of
$\M_-$ in the case of cubic fourfolds and the rank of $\LLL_+$ in
the case of K3-surfaces. Respectively, as it follows from the
Lefschetz trace formula, $r$ is equal to
$\frac12\big(\chi(X)-\chi(X(\R))-4\big)=
11+\frac12\big(1-\chi(X(\R))\big)$ in the first case, and
$\frac12\big(\chi(Y) +\chi(Y(\R))-4\big)= 10+\frac12\chi(Y(\R))$
in the second.

In $\Gamma_{K4}$ as well as in $\Gamma_{K3}$ there are 10 pairs of
vertices-twins having the same coordinates $(r,d)$, and on Figures
1 and 2 we place one of the vertices in each pair slightly above
the other vertex (one may consider this as introducing a third
coordinate).

The vertices in $\Gamma_{K4}$ and in $\Gamma_{K3}$ which are
marked by $\bullet$ have lattices of the type I. The other
vertices are marked by $\circ$ and have lattices of the type II
(we recall the definition of the types in \ref{mark}). Note that
in a pair of vertices-twins with the same coordinates $(r,d)$ one
vertex belongs to the type I and the other to type II. Among the
vertices which are completely determined by $(r,d)$ the majority
belong to the type II, but the five bottom ``local minima''
vertices belong to the type I.

It will be convenient for us  to introduce, instead of $(r,d)$, a
new coordinate system $(i,j)$ centered at the top vertex, namely,
the system where $i=\frac12(22-r-d)$ is the number of down-left
and $j=\frac12(r-d)$ of down-right moves along edges needed to
reach a given vertex from the top. We will call such down-left
moves {\it L-moves} and the down-right moves {\it R-moves}. The
moves in the opposite directions will be called respectively {\it
$L^{-1}$-moves} and {\it $R^{-1}$-moves}. A vertex with
coordinates $(i,j)$ will be denoted $C^{i,j}$ (in the case of the
K3-graph as well as for the K4-graph) if it is the only vertex
with such coordinates. In the case of two vertices with
coordinates $(i,j)$ we use notation $C^{i,j}_I$ for the vertex of
type I and $C^{i,j}$ for the other vertex. The set formed by the
vertices $C^{i,j}$ will be called {\it the principal series} of
vertices, and the vertices $C^{i,j}_I$ will be called the {\it
special vertices}.
 We say that a vertex
of the K3-graph and a vertex of the K4-graph {\it correspond to
each other} if they are of the same type and have the same $(i,j)$
coordinates (i.e., the both have the same notation, either
$C^{i,j}$ or $C^{i,j}_I$).

The walls (the facets) between the deformation chambers are
divided into {\it L-walls} and {\it R-walls} ({\it L-facets} and
{\it R-facets}), depending on the edges they represent in
$\Gamma_{K4}$.

\subsection{Immediate topological consequences of deformation
classification}\label{immediate} It is easy to construct in any
dimension $n$ a non-singular real cubic hypersurface $X\subset
P^{n+1}$ whose real locus $X(\R)$ is diffeomorphic to $\Rp{n}$, as
well as a cubic with $X(\R)=\Rp{n}\+S^n$. Namely, it is sufficient
to perturb a singular cubic which splits into a hyperplane and a
quadric; if the quadric has an empty real locus, then we obtain
$X(\R)=\Rp{n}$, and if it has a real locus disjoint from the
hyperplane, then $X(\R)=\Rp{n}\+S^n$. Clearly, in each of these
two cases the models constructed in such a way are deformation
equivalent to each other.

\lemma\label{disconnected} If the real locus $X(\R)$ of a
non-singular real cubic $n$-fold $X$, $n\ge1$, is disconnected
then $X(\R)$ is diffeomorphic to a disjoint union of $\Rp{n}$ and
$S^n$. Such cubics form a single deformation component.
\endlemma

\demo{Proof} As it follows from B\'ezout theorem, $X(\R)$ realizes
a non-trivial class in the homology group $H_n (\Rp{n+1};\Z/2)$.
So does also one of its connected components, $C_1$. Any other
connected component, $C_2$, should realize a trivial class since
$C_1\cap C_2=\oo$. The complement $\Rp{n+1}\smallsetminus C_2$
splits in two connected components one of which contains $C_1$,
and we choose a point $p$ in the other component. Any real line
through $p$ intersects $C_1$ at least once and $C_2$ at least
twice. So (again by B\'ezout theorem) $X(\R)$ cannot have other
components and the central projection $X(\R)\to\Rp{n}$ from $p$
gives a 1-1 map $C_1\to\Rp{n}$ and a double covering
$C_2\to\Rp{n}$. This proves the first statement of the Lemma.

The same arguments show that if a real cubic $(n-1)$-fold (in a
real projective $n$-space) is disconnected and has a real singular
point, then its real locus is the union of this point and a real
projective hyperplane disjoint from it. Therefore, at any point of
$C_2$ the hyperplane tangent to $C_2$ intersects $C_2$ exclusively
by this tangency point. Let us shift
 the
tangency hyperplane to make it a real hyperplane $H$ disjoint from
$C_2$, pick up a small ellipsoid $E$ inside $C_2$ and consider the
cubic $X'$ which is a small real perturbation of $H\cup E$. Such a
cubic $X'$ is deformation equivalent to $X$, since they can be
joined by a straight deformation path, $tf+(1-t)g$, where the
degree $3$ defining equations $f=0$ and $g=0$ of $X$ and $X'$,
respectively, are chosen in such a way that $fg>0$ at the solid
ellipsoid bounded by $E$ in $P^{n+1}(\R)$: under this sign
convention the cubics $tf+(1-t)g=0, t\in [0,1],$ remain
disconnected and, hence, non-singular. Now, the second statement
of the Lemma follows from the mutual deformation equivalence of
all the model disconnected cubics.\qed\enddemo

\corollary\label{basecase} {\rm (a)} If a non-singular real cubic
fourfold has real part $X(\R)=\Rp{4}$, then it represents the
vertex $C^{0,0}$. {\rm (b)} If a non-singular real cubic fourfold
has real part $X(\R)=\Rp{4}\cup S^4$, then it represents the
vertex $C^{1,0}_I$.
\endcorollary

\demo{Proof} In the case {\rm (a)}, the relations $r=
11+\frac12\big(1-\chi(X(\R))\big)=11$, $d=\frac12\big(
27-b_*(X(\R))\big)=11$ give $i=\frac12(22-11-11)=0$,
$j=\frac12(11-11)=0$, and we see from the graph $\Gamma_{K4}$ that
$C^{0,0}$ is the unique suitable vertex. Similarly, in the case
{\rm (b)} we obtain $i=1$, $j=0$ and then notice in addition that
the only connected four-manifold that can be obtained from
$\Rp{4}\cup S^4$ by a single Morse modification is $\Rp{4}$, while
as is seen on the graph $\Gamma_{K4}$ the vertex $C^{1,0}$ is
adjacent to two vertices. Thus, the corresponding vertex is
$C^{1,0}_I$. \qed\enddemo

\subsection{Markings and eigenlattices (\cite{FK1,FK2})}\label{mark} By the
lattice $\M(X)$ of a non-singular cubic fourfold $X$ we mean the
middle cohomology group $H^4(X;\Z)$ endowed with the intersection
form. This is an odd unimodular lattice of signature $(21,2)$. The
{\it polarization class} $h_X\in \M(X)$ (which is realized by the
intersection of $X$ with a projective $3$-space) has square
$h_X^2=3$ and is {\it characteristic} which means as usual that
$xh_X=x^2\mod2$ for all $x\in \M(X)$. This implies that there
exists a lattice isomorphism between $\M(X)$ and $\M=3I+2U +2E_8$
such that $h_X$ is mapped into $h=(1,1,1)\in3I$; in particular,
the primitive sublattice $\M^0=\{x\in \M\,|\, xh=0\}$ becomes
identified with $A_2+2U +2E_8$. Such a choice of a lattice
isomorphism $\phi\:\M(X)\to \M$ will be called {\it a marking of
$X$}, and a pair $(X,\phi)$ will be called {\it a marked cubic
fourfold}.

If a cubic $X$ is defined over reals, then the complex conjugation
induces a lattice involution $c_X\:\M(X)\to \M(X)$ such that
$c_X(h_X)=h_X$. This involution gives rise to {\it the
eigenlattices} $M_\pm(X)=\{x\in \M(X)\,|\,c_X(x)=\pm x\}$ and
$\M_\pm^0(X)=\{x\in \M^0(X)\,|\,c_X(x)=\pm x\}$. Here,
$\M^0_-(X)=\M_-(X)$ since $h_X\in \M_+(X)$. A marking
$\phi\:\M(X)\to\ \M$ transforms the involution $c_X$ and its
eigenlattices $\M_\pm(X)$, $\M_\pm^0(X)$ into an involution on
$\M$ and the eigenlattices of this induced involution.

An involution $c:\M\to \M$ is induced by $c_X$ for some marked
real cubic fourfold $(X,\phi)$ if and only if the involution $c$
is {\it geometric}, where the latter means that $c(h)=h$ and each
of the eigenlattices $\M_\pm=\{x\in \M\,|\,c(x)=\pm x\}$ and
$\M_\pm^0=\{x\in \M\,|\,c(x)=\pm x, \, xh=0\}$ is of negative
inertia index equal to $1$. In fact, the isomorphism class of the
pair $(\M(X),c_X)$ is an invariant which allows to distinguish the
coarse deformation classes of real cubic fourfolds $X$. Moreover,
the isomorphism class of $(\M(X),c_X)$ is determined by the
isomorphism class of any of the eigenlattices $\M^0_+(X)$ and
$\M_-(X)$.

The isomorphism classes of lattices $\M^0_+$ and $\M_-$
corresponding to geometric involutions are listed in the following
tables 1---3.

\midinsert \line{\vtop{\hsize 5.2cm\hskip-4mm
$\matrix\text{Table 1.}\\ \text{$\M_+^0$ for cubic fourfolds of type $C^{i,j}$}\\
\boxed{\matrix
 \format\l&\quad\r&\l\\
C^{0,j}& 
-A_1+(9-j) A_1&+\la6\ra\\
 C^{1,j}& 
-A_1+(9-j) A_1&+A_2\\
C^{2,j}& 
U+(9-j) A_1&+A_2\\
 C^{3,j}& 
U+(6-j) A_1&+A_2+D_4\\
 C^{4,j}& 
-A_1+(5-j) A_1&+\la6\ra+E_8\\
C^{5,j}& 
-A_1+(5-j) A_1&+A_2+E_8\\
 C^{6,j}& 
U+(5-j) A_1&+A_2+E_8\\
 C^{7,j}& 
U+(2-j) A_1&+A_2+D_4+E_8\\
 C^{8,j}& 
-A_1+(1-j) A_1&+\la6\ra+2E_8\\
 C^{9,j}& 
-A_1+(1-j) A_1&+A_2+2E_8\\
 C^{10,j}& 
U+(1-j) A_1&+A_2+2E_8\\
\endmatrix
}\endmatrix$}
 \kern2pt\vtop{\hsize 6.4cm
$\matrix\text{Table 2.}\\ \text{$\M_-$ for cubic fourfolds of type
$C^{i,j}$}\\ \boxed{\matrix
  \format\l&\quad\r&\l\\
 C^{i,0} &-A_1+(10-i)A_1\\
 C^{i,1} &U+(10-i)A_1\\
 C^{i,2} &U+(7-i)A_1&+D_4\\
 C^{i,3} &-A_1+(6-i)A_1&+E_7\\
 C^{i,4} &-A_1+(6-i)A_1&+E_8\\
 C^{i,5} &U+(6-i)A_1&+E_8\\
 C^{i,6} &U+(3-i)A_1&+D_4+E_8\\
 C^{i,7} &-A_1+(2-i)A_1&+E_7+E_8\\
 C^{i,8} &-A_1+(2-i)A_1&+2E_8\\
 C^{i,9} &U+(2-i)A_1&+2E_8\\
\endmatrix} \\ \\ \endmatrix$}}
\endinsert

\midinsert \centerline{\hskip10mm $\matrix\text{Table 3.}\\
\text{$\M_+^0$ and $\M_-$ for cubic
fourfolds of type $C^{i,j}_I$}\\
\boxed{\matrix
C^{i,j}_I&\M_+^0&\M_-\\
\text{-----}&\text{-----------------------}&\text{------------------}\\
 C^{0,3}_I&U(2)+E_6(2)&U(2)+3D_4\\
 C^{1,8}_I &U(2)+A_2&U(2)+2E_8\\
 C^{1,4}_I &U+E_6(2)&U+3D_4\\
 C^{2,5}_I &U(2)+A_2+D_4&U(2)+D_4+E_8\\
 C^{3,2}_I &U(2)+A_2+2D_4&U(2)+2D_4\\
 C^{4,3}_I &U+A_2+2D_4&U+2D_4\\
 C^{5,4}_I &U(2)+A_2+E_8&U(2)+E_8\\
 C^{6,1}_I &U(2)+A_2+D_4+E_8&U(2)+D_4\\
 C^{9,0}_I &U(2)+A_2+2E_8&U(2)\\
 C^{2,1}_I &U+A_2+E_8(2)&U+E_8(2)\\
 C^{1,0}_I &U(2)+A_2+E_8(2)&U(2)+E_8(2)\\
\endmatrix}\endmatrix$}
\endinsert

We say that a real cubic $X$ is {\it of type} I if the involution
$c_X \:\M(X)\to \M(X)$ is {\it even}, that is if $x\cdot
c_X(x)+x^2=0\mod2\Z$ for any $x\in \M(X)$. As is known, it happens
if and only if $X(\R)$ realizes the class $h_X$ in the middle
$(\Z/2)$-homology of $X$.

The above definition of type I is not specific to cubics, but is
applied to any non-singular real algebraic variety. In particular,
if $Y$ is a $K3$-surface, it is of type I if and only if  $Y(\R)$
realizes zero in the middle $(\Z/2)$-homology of $Y$.

\subsection{Central projection K3-K4 correspondence \cite{FK1}}\label{central-correspondence}
The edges of the K4-graph can be interpreted as coarse deformation
classes of 6-polarized K3-surfaces via the following {\it central
projection correspondence}. Let $X_0$ denote a nodal real cubic
fourfold representing an edge. In an affine chart centered at the
node, the cubic $X_0$ is defined as $\{f_2+f_3=0\}\subset\R^5$
where $f_2$ and $f_3$ are some quadratic and, respectively, cubic
homogeneous polynomials in $x_0,\dots,x_4$. Non-degeneracy of the
node means that the quadric $\{f_2=0\}\subset\Rp4$ is non-singular
and absence of other singularities in $X_0$ is equivalent to
transversality  of the intersection $\{f_2=f_3=0\}\subset\Rp4$,
which is therefore a non-singular real 6-polarized K3-surface.

\theorem\label{k3-k4correspondence} \roster\item A vertex $v_{K3}$
of the K3-graph corresponds to a vertex $v_{K4}$ of the K4-graph
if and only if the eigenlattice $\M_-(X)$ of a cubic fourfold $X$
representing $v_{K4}$ is isomorphic to the eigenlattice
$-\LLL_+(Y)$ of a K3-surface $Y$ representing $v_{K3}$. \item The
edges of the K4-graph are in one-to-one correspondence with the
coarse deformation classes of 6-polarized K3 surfaces $Y$. Namely,
edges with the leftmost vertex \ $v_{K4}$ correspond to the
classes of $(Y,l)$, such that $Y$ represent $v_{K3}$ and
$l\in\LLL_-(Y)$, $l^2=6$.\endroster
\endtheorem

\subsection{The period map (\cite{Laza}, \cite{Looijenga}, \cite{Voisin})}
The non-zero Hodge numbers in dimension four for any non-singular
cubic fourfold $X\subset P^5$ are $h^{3,1}=h^{1,3}=1$ and
$h^{2,2}=21$. Given a marking $\phi: (\M(X), h_X)\to (\M,h)$, the
complex line $\phi(H^{3,1}(X))\subset \M^0\otimes\C$ is isotropic
and has negative pairing with the conjugate (and thus, also
isotropic) line $\phi(H^{1,3}(X))=\overline{\phi(H^{3,1}(X))}$,
that is to say, $w^2=0$, and $w\overline w<0$, (and thus
$\overline w^2=0$) for all $w\in \phi(H^{3,1}(X))$. Writing
$w=u+iv$, $u,v\in \M^0\otimes\R$, we can reformulate it as
$u^2=v^2<0$ and $uv=0$, which implies that the real plane $\la
u,v\ra\subset \M^0\otimes\R$ spanned by $u$ and $v$ is negative
definite and bears a natural orientation given by $u=\Re w, v=\Im
w$. Note that the orientation determined similarly by the complex
line $\phi(H^{1,3}(X))\subset \M^0\otimes\C$ is the opposite one.

The line $\phi(H^{3,1}(X))\subset \M^0\otimes\C$ specifies a point
$\Omega(X)\in P(\M^0\otimes\C)$ (as usual, $P$ states for the
projectivization) called the {\it period point of $(X,\phi)$}.
This period point belongs to the quadric $\QQQ= \{w^2=0\}\subset
P(\M^0\otimes\C)$, and more precisely, to its open subset,
$\widehat{\DD}=\{w\in \QQQ\,|\,w\overline w<0\}$.
 This subset has two connected components,
which are exchanged by the complex conjugation (this reflects also
switching from the given complex structure on $X$ to the complex
conjugate one).

The orthogonal projection of a negative definite real plane in
$\M^0\otimes\R$ to another such one is non-degenerate. Thus, to
select one of the two connected components of $\widehat{\DD}$ we
fix an orientation of negative definite real planes in
$\M^0\otimes\R$ which is constrained to be preserved by the
orthogonal projection. We call it the {\it prescribed orientation}
and restrict the choice of markings to those for which the
orientation of $\phi(H^{3,1}(X))$ defined by the pairs $u=\Re w,
v=\Im w$ for $w\in\phi(H^{3,1}(X))$ is the prescribed one. We
denote this component by $\DD$ and call it the {\it period
domain}. By $\Aut^+(\M^0)$ we denote the group of those
automorphisms of $\M^0$ which preserve the prescribed orientation
(and thus preserve $\DD$). We put
$\Aut^-(\M^0)=\Aut(\M^0)\setminus\Aut^+(\M^0)$. This complementary
coset consists of automorphisms exchanging the connected
components of $\widehat{\DD}$.

On the other hand we have the projective space $P_{4,3}$ formed by
all cubic fourfolds, which splits into the {\it discriminant
hypersurface} $\D_{4,3}$ formed by singular cubics and its
complement, $\CC=P_{4,3}\setminus \D_{4,3}$. Let $\CC^\sharp$
denote the space of marked non-singular cubics. The natural
projection $\CC^\sharp\to\CC$ is obviously a Galois covering with
the deck transformation group $\Aut^+(\M^0)$. The above
conventions define the {\it period map} $\per\: \CC^\sharp\to
\DD$, $(X,\phi)\mapsto \phi(H^{3,1}(X))$.

The above definitions extend naturally to cubic fourfolds with
simple singularities. In particular, the lattice $H^4(X)$ of a
cubic fourfold with simple singularities is torsion free and
admits an isometric embedding $H^4(X)\to\M$ whose orthogonal
complement is isometric to $\bigoplus_i \M_{x_i}(X)$, the sum of
the Milnor lattices $\M_{x_i}$ over all singular points $x_i$ of
$X$. On the other hand, $H^4(X)\otimes\C$  carries a pure Hodge
structure with $h^{3,1}=h^{1,3}=1$. By Riemann extension theorem
and the finiteness of the monodromy groups of simple
singularities, the Galois covering $\CC^\sharp\to\CC$ extends (in
an unique way) to a ramified Galois covering $\CC_s^\sharp\to
\CC_s$ where $\CC_s\subset P_{4,3}$ is the space of cubic
fourfolds with simple singularities. The covering space
$\CC_s^\sharp$ is non-singular. It is the space of {\it marked
cubic fourfolds with simple singularities}, where by a marking of
$X\in \CC_s$ we understand, in accordance with the previous
definitions, a respecting the prescribed orientation (of negative
definite planes) isometric embedding $(H^4(X), h_X)\to(\M, h)$
whose orthogonal complement is isometric to the sum of the Milnor
lattices over all singular points of $X$. The extended period map
$\per\: \CC_s^\sharp\to \DD$, $(X,\phi)\mapsto \phi(H^{3,1}(X))$,
is holomophic due to the corresponding Griffiths theorem.

Consider, on the period space side, the reflection $R_v$ in $
\M^0\otimes\C$ across the mirror-hyperplane $H_v=\{x\in
\M^0\otimes\C\, |\, xv=0\}$ defined as $x\mapsto
x-2\frac{xv}{v^2}v$, and note that it preserves the lattice $\M^0$
invariant if $v\in \M^0$ is such that $v^2=2$, or such that
$v^2=6$ and $xv$ is divisible by $3$ for all $x\in \M^0$. We call
these two types of lattice elements {\it $2$-roots} and {\it
$6$-roots} respectively, and denote their sets by $V_2$ and $V_6$.
Note that $R_v\in\Aut^+(\M^0)$ for any $v\in V_2\cup V_6$. If
$v\in V_2$, then the reflection $R_v$ extends (as a reflection) to
$\M$ and $h$ is preserved by this extension. By contrary if $v\in
V_6$, the reflection $R_v$ does not extend to a reflection in
$\M$, and moreover, the unique extension of $R_v$ to $\M$ maps $h$
to $-h$. On the other hand, if $v\in V_6$ then the
\emph{anti-reflection} $-R_v$ extends to an isometry of $\M$
preserving $h$. This extension is the anti-reflection with respect
to the 2-plane generated by $h$ and $v$. In particular, it
represent also an element of $\Aut^+(\M^0)$.

The union of the mirrors $H_v$ for all $v\in V_2$ gives after
projectivization a union $\Cal H_\D\subset P(\M^0\otimes\C)$ of
hyperplanes, and a similar union of $H_v$ for all $v\in V_6$ gives
another union of hyperplanes, $\Cal H_\infty \subset
P(\M^0\otimes\C)$.

\theorem["Surjectivity" of the period map]\label{C-surj} The image
of the period map $\per : \CC_s^\sharp\to \DD$ is the complement
of $\Cal H_\infty$, and the fibers of $\per$ are $PGL(6,\C)$
orbits. For a given $p\in \DD\setminus \Cal H_\infty$, the 2-roots
$\delta\in V_2$ such that $H_\delta$ contains $p$ form an elliptic
root system, whose irreducible components are of types $A, D$, and
$E$. These components generate the Milnor lattices of the singular
points of a cubic with the period $p$. \qed\endtheorem

\remark{{\eightit Remark}}{\eightrm The description of the fibers
given in the first part of the statement is equivalent to what is
called usually the {\eightit injectivity of the period map} (here,
at the level of cubics with simple singularities). The second part
of the statement includes an intermediate, in some sense,
statement that the variations of a singular cubic fourfold contain
a simultaneous versal deformation of the singularities if all of
them are simple.}
\endremark\vskip0.05in

\subsection{The period map in a real setting (\cite{FK2})}
Let us fix a geometric involution $c:\M\to \M$, see \ref{mark}. A
{\it real $c$-marked nonsingular cubic fourfold {\rm
(respectively,} cubic fourfold with simple singularities{\rm )}}
is, by definition, a real non-singular cubic fourfold
(respectively, cubic fourfold with simple singularities) equipped
with a marking $\phi$ such that $\phi\circ\conj^*=c\circ\phi$. If
such a $c$-marking exists, the cubic fourfold is said to be of
{\it homological type} $c$.

We denote by $\CC^c_\R\subset\CC_\R$ (respectively,
$\CC^c_{s,\R}\subset\CC_{s,\R}$) the set of real nonsingular cubic
fourfolds (respectively, cubic fourfolds with simple
singularities) of homological type $c$, and by $\CC_\R^{c\sharp},
\CC^{c\sharp}_{s,\R}$ the respective sets of $c$-marked real cubic
fourfolds. The two latter sets can be seen as the real parts of
$\CC^\sharp$ and $\CC^\sharp_s$ with respect to the involution
which sends $(X,\phi)\in\CC^\sharp$ to
$(\conj(X),c\circ\phi\circ\conj^*)$.

Let us extend $c$ to a complex linear involution on $\M\otimes \C$
and denote also by $c$ the induced involutions on $
\M^0\otimes\C$, $P=P( \M^0\otimes\C)$, and $\widehat{\DD}$. Note
that $c$ permutes the two components $\DD$ and $\overline{\DD}$ of
$\widehat{\DD}$, and thus, $\overline c(\DD)=\DD$, where
$\overline c\: \M^0\otimes\C\to\M^0\otimes\C$ is the composition
of $c$ with the complex conjugation in $ \M^0\otimes\C$.

Let $\widehat{\DD}_\R^{c}$ and $\DD_\R^{c}$ denote the fixed point
set of $\overline c$ restricted to $\widehat{\DD}$ and $\DD$. The
second set, $\DD_\R^{c}$, consists of the lines generated by
$w=u_++iu_-$ such that $u_\pm\in \M^0_\pm(c)\otimes\R$,
$u_+^2=u_-^2<0$, and the orientation $u_+,u_-$ is the prescribed
one. Since $c$ is geometric, both ${\DD}_\R^{c}$ and its (trivial)
double covering $\widehat{\DD_\R^{c}}$ are nonempty.

As it follows from definitions, the period point of a c-marked
real cubic fourfold belongs to
$\DD_\R^{c}=\{x\in\DD\,|\,c(x)=\overline x\}$. Therefore, we may
speak of a {\it real period map} $\per_\R :
\CC_{s,\R}^{c\sharp}\to \DD_\R^{c}$ and call $\DD_\R^{c}$ the {\it
real period domain of real $c$-marked cubic fourfolds}.

\theorem["Surjectivity" of the period map in a real
setting]\label{rsurjectivity} The image of the real period map
$\per_\R : \CC_{s,\R}^{c\sharp}\to \DD_\R^{c}$ is the complement
of $\Cal H_\infty\cap \DD_\R^{c}$, and the fibers of $\per_\R$ are
$PGL(6,\R)$ orbits. For a given $p\in \DD_\R^{c} \setminus \Cal
H_\infty$, the 2-roots $\delta$ such that $H_\delta$ contains $p$
form an elliptic root system invariant under $c$-action. Its
irreducible components which are invariant under $c$-action
generate the Milnor lattices of the  real singular points of a
cubic with the period $p$.
\endtheorem

As usual, such a real statement can be deduced from the
corresponding complex statement, which means here from
Theorem~\ref{C-surj}. To do it we use the approach applied in
~\cite{ACT} in a similar situation.

\demo{Proof}
 The action of $PGL(6,\C)$ on $\CC_s^\sharp$ and the map
$\CC_s^\sharp/PGL(6,\C)\to \DD\setminus\Cal H_\infty$ induced by
the period map are proper. Therefore, there is a complete
Riemannian metric on $\CC_s^\sharp$ which is invariant under the
action of $PGL(6,\C)$ and descends to a complete Riemannian metric
on $\DD\setminus\Cal H_\infty$. Taking the average, one can make
these metrics invariant under the complex conjugation. Then, the
induced Riemannian metric on $\CC_{s,\R}^{c\sharp}$ is complete
and invariant under the action of $PGL(6,\R)$. Since, in addition,
the action of $PGL(6,\C)$ has finite stabilizers, the map $
\CC_{s,\R}^{c\sharp}/PGL(6,\R)\to \DD^c_\R$ (induced by the period
map) is proper, which (together with the local Torelli theorem)
implies its surjectivity. Its injectivity follows from the local
Torelli and the injectivity on the subspace of (marked real)
nonsingular cubics.

The second statement is a straightforward consequence of the
second statement of Theorem~\ref{C-surj}. \qed\enddemo

\section{Arithmetics and topology under wall crossings}\label{AiT}

\subsection{Indices of the Morse modifications under the
facet crossing} Consider non-singular real cubic fourfolds $X$ and
$X'$ representing adjacent deformation components $C$ and $C'$. As
we connect $X$ with $X'$ by a continuous path which crosses (once)
some facet, $\Cal F$, the real locus $X(\R)$ experiences a Morse
modification of index $0\le q\le5$. If we follow this path in the
opposite direction, then we observe that
 $X(\R)$ is obtained from $X'(\R)$ by a Morse
modification of index $p=5-q$. In particular,
$\chi(X'(\R))-\chi(X(\R))$ is $2$ if $q$ is even and $-2$ if odd.
We have also $d(X')-d(X)=\pm1$ (see Figure 1 or subsection
\ref{w-crossing-arithms} below), where $d$ as before denote the
discriminant rank of the eigenlattice $\M_-$. If $d(X')<d(X)$,
then we define the {\it index of facet $\Cal F$} as $\ind(\Cal
F)=q$; in the case $d(X')>d(X)$ we obtain respectively $\ind(\Cal
F)=p=5-q$.

\lemma\label{LRparity} \roster\item L-facets have even index and
R-facets have odd index. \item The core sphere of the Morse
modification for L-moves and for R-moves are null-homologous in
$H_{q-1}(X(\R);\Z/2)$.\item The core spheres of the Morse
modification for $L^{-1}$-moves and for $R^{-1}$-moves are
homologically non-trivial in $H_{p-1}(X(\R);\Z/2)$.
\endroster
\endlemma

\demo{Proof} L-moves increase and R-moves decrease the Euler
characteristic $\chi(X(\R))$, since $r=\rank(M_-)$ plays the role
of the horizontal coordinate on Figure 1 (see \ref{K4}). This
implies (1). After L-moves and R-moves the discriminant rank
(i.e., the vertical coordinate)
$d(X)=\frac12\big(b_*(X)-b_*(X(\R))\big)$ is decreasing and thus
$b_*(X(\R))$ is increasing. Conversely, after $L^{-1}$-moves and
$R^{-1}$-moves $b_*(X(\R))$ is decreasing. This implies (2) and
(3). \qed\enddemo

\lemma\label{LRindex} \roster\item An R-facet has either index
$1$, or $3$; \item the L-facets between $C^{0,0}$ and $C^{1,0}_I$
have either index $0$, or $4$; \item all other L-facets have index
$2$.
\endroster\endlemma

\demo{Proof} R-facets cannot have index 5, as follows from the the
K4-graph and connectedness of $X(\R)$ for $X\notin C^{1,0}_I$.
This implies (1).
 A spherical component of $x\in C^{1,0}_I$ may be only obtained by
 a Morse modification of index $0$ or $4$, which implies (2).

To prove (3), we observe that if L-facet has index $4$, then by
\ref{LRparity}(2) its core sphere should be homologically trivial
($\roman{mod}\ 2$) and thus the Morse modification increases the
number of components of $X(\R)$, like in the case of an index $0$
modification. \qed\enddemo

\subsection{Cuspidal cubics and their
perturbations}\label{cusp-perturb} In this section we examine
cubics of arbitrary dimension $n$.

Assume that a wall separating two deformation components contains
two facets adjacent through a cuspidal stratum. Consider a generic
point, $\a\in\D(\R)$, of this stratum, which represents a real
cubic $X_\a$ with a cusp and no other singular points. Near the
cusp, in a suitable affine chart, the equation of $X_\a$ is
$z^3+f_2+zg_2+f_3=0$, where $f_2=f_2(x,y)=\sum_{i=1}^p
x_i^2-\sum_{j=1}^q y_j^2$, $p+q=n$, and $g_2=g_2(x,y),
f_3=f_3(x,y)$ are some homogeneous polynomial of degree 2 and 3,
respectively. Let us include $X_\alpha$ in a two-dimensional
linear system of cubics $X_{b,c}=\{f_{b,c}=0\}$, $f_{b,c}=z^3+ bz
+c+f_2+zg_2+f_3$, $b,c\in\R$. This linear system yields a versal
deformation of the cusp singularity, moreover, there exists a
local change of $(x,y)$-coordinates reducing, without change of
$z$-coordinate, the polynomials $f_{b,c}$ to $z^3+ bz +c+f_2$. It
implies that near $\alpha$ our linear system intersects with the
discriminant $\D$ transversally along a cuspidal curve
$\D_{b,c}=\{4b^3+27c^2=0\}$. In particular, near $\a$ the
hypersurface $\D(\R)$ is a topological manifold which is split by
the cuspidal stratum (corresponding to $b=c=0$) into a union of
two adjacent facets: one corresponding to $c>0$ and the other to
$c<0$.

To treat below the cubics of any dimension we have to extend the
notion of index of facets in accordance with our definition for
$n=4$. We do it via the following "coorientation" convention: the
total $\Z/2$-Betti number $b_*(X_\a(\R))$ of the real locus of the
cubic representing a generic point, $\a\in P_{n,3}(\R)$, in the
space of real cubics should increase as $\a$ crosses a facet in
the positive direction.

This convention does coorient the facets adjacent to the cuspidal
strata as it is shown on Figure 3: the normal vectors looks into
the ``thinner'' region bounded by the cusp-shaped discriminant.
This is a part of Lemma \ref{through-cusp} below, where we show
that such crossing a facet adds a handle to $X_\a$. For this aim,
we look at the part $B_{b,c}=X_{b,c}(\R)\cap B_\e$ of $X_{b,c}$
inside a Milnor ball $B_\e$ (of radius $\e$), where
$|b|,|c|<\!<\e<\!<1$.

\lemma\label{through-cusp} \roster \item The facet $c>0$ of the
discriminant $\D_{b,c}=\{4b^3+27c^2=0\}$,  has index $q$, and its
facet $c<0$ has index $p$. \item If $4b^3+27c^2<0$, then $B_{b,c}$
is a smooth $n$-ball properly embedded in the $(n+1)$-ball $B_\e$.
\item If $4(b')^3+27(c')^2>0$, then $B_{b',c'}$ is isotopic to
$B_{b,c}\#(S^p\times S^q)$, i.e., ambient connected sum of the
above $B_{b,c}$ with an unknotted handle.
\endroster
\endlemma

\midinsert \topcaption{Figure 3. Indices of facets adjacent
through a cuspidal stratum}
\endcaption\epsfbox{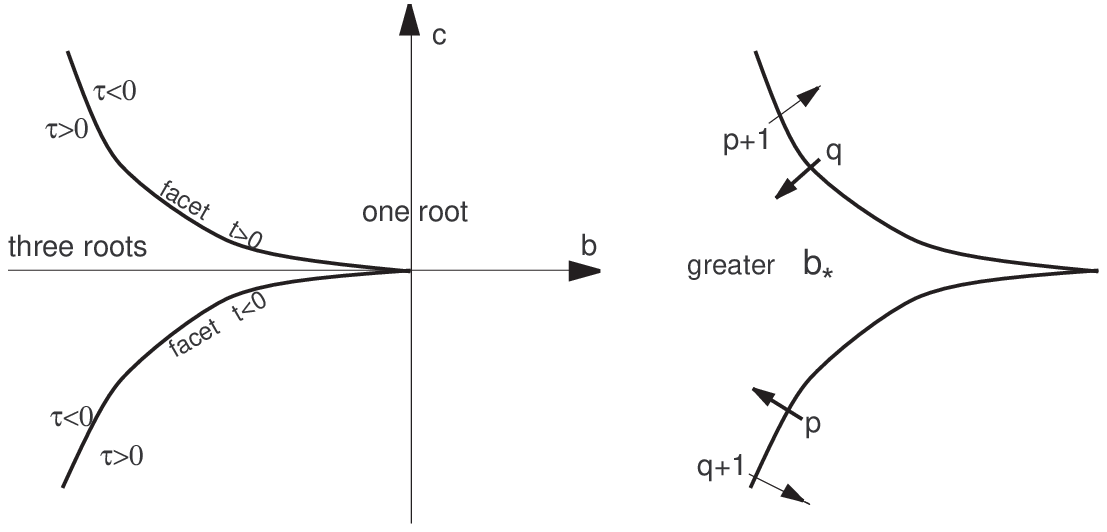}
\botcaption{\eightpoint \noindent The base case is on the
left-hand side figure. Its stabilization is on right-hand side
figure, where the arrows show the direction of crossing the
facets, and next to the arrow stand the indices of the
corresponding Morse modifications.}\endcaption
\endinsert

\proof Let us check the base of induction, $n=0$, that is the case
of one variable polynomial $z^3$ (in this case, $p=q=0$). It is
sufficient to examine the behavior of $z^3-3t^2z + 2t^3-\tau$. For
$\tau=0$ this polynomial has a double real root $z=t$ and a simple
real root $z=-2t$. If $t>0$, then this double root disappears for
$\tau<0$ and turns into two real roots for $\tau>0$, hence, the
index of this facet is equal to $0$ ($=q$). If $t<0$, then the
double root disappears for $\tau>0$ and turns into two real roots
for $\tau<0$; therefore, to get an increase of the total Betti
number (here it is the number of real roots) with respect to
growing $\tau$ we alternate the sign of the equation, and thus
find that the index of this facet is equal to $0$ ($=p$).

Since the stabilization of the equation (i.e., adding or
subtracting the square of a new variable) does not change the
total Betti number of its locus in the Milnor ball, it is now
straightforward to complete the proof of (1) by induction.

Parts (2) and (3) are also proved by induction. Namely, consider
regions $B_{b,c}^\pm=\{(x,y,z)\in B_\e\,|\,\pm f_{b,c}\ge0\}$
separated by $B_{b,c}$. After we increase $p$ or $q$, we get a new
locus $B_{b,c}$, which is obtained by a suspension, that is can be
identified with the ``double'' of the corresponding region
$B_{b,c}^\pm$ (the ``double'' in this context is obtained from two
copies of $B_{b,c}^\pm$ glued along $B_{b,c}$). This immediately
implies (2). For proving (3) we should additionally notice that a
handle in the ambient connected sum $B_{b',c'}=B_{b,c}\#(S^p\times
S^q)$ is not knotted provided that the spheres $S^p\times\pt$ and
$\pt\times S^q$ (which are the vanishing cycles corresponding to
the facets $c>0$ and $c<0$ respectively) bound disjoint balls
(each in the corresponding region $B_{b',c'}^\pm$). One can easily
see that the suspension (which is increasing $p$ or $q$) preserves
this property of the spheres.
\endproof

Returning to dimension $n=4$ we get the following result.

\corollary\label{cusp-stratum}\roster \item In a pair of
cuspidally adjacent $R$-facets one has always index 1 and another
index 3. \item  In a pair of cuspidally adjacent $L$-facets
separating $C^{0,0}$ and $C^{1,0}_I$ one has index $0$ and another
index $4$. \qed\endroster
\endcorollary

\subsection{Arithmetics of the wall-crossing (\cite{FK1, FK2})}\label{w-crossing-arithms}
Before sticking to the case of cubic fourfolds, let us make a
somewhat general remark concerning real non-singular projective
hypersurfaces $X$ of any dimension $n$, and their involutions
$c\:H^n(X)\to H^n(X)$ induced by the complex conjugation. It is a
well known (and simple) fact that adjacency of the deformation
components of a pair of such hypersurfaces, $X_\pm\in
P_{n,d}(\R)\setminus\Delta(\R)$, implies that their involutions,
$c_\pm$, are related by the Picard-Lefschetz transformation.
Namely, we can connect $X_\pm$ by a family $X_t$ of non-singular
hypersurfaces, $X_t\in P_{n,d}\setminus\Delta$, which follows a
path in the real parameter space $P_{n,d}(\R)$ with exception of
an arc, $\gamma$, in the complex domain, $P_{n,d}\setminus
P_{n,d}(\R)$. Then the composition $c_-\circ c_+$ turns out to be
the monodromy along the loop $\bar\gamma^{-1} *\gamma$. If
$\gamma$ is a small arc which turns around the obstructing facet
$\Cal F$ counterclockwise this monodromy is, by definition, the
Picard-Lefschetz transformation.

In even dimensions $n$, the Picard-Lefschetz transformation maps
$x\in H^n(X_+)$ to $R_v(x)=x- (-1)^\frac{n}2 (v\cdot x)v$ where
$v$ is the vanishing (co)-cycle. Therefore, the two involutions
$c_\pm$ coincide in the orthogonal complement of $v$, while $v$
jumps from $\pm1$-eigenspace of $c_+$ to the opposite
$\mp1$-eigenspace of $c_-$. Thus, by the Lefschetz trace formula,
$\chi(X_+(\R))-\chi(X_-(\R))$ equals $2$ if $v$ belongs to the
$(+1)$-eigenspace, and $-2$ if $v$ belongs to the
$(-1)$-eigenspace.
 It is also well-known (and straightforward) that
the discriminants, $d(X_+)$ and $d(X_-)$, of $X_\pm$ should differ
by $1$. Namely, $d(c_+)-d(c_-)=1$ if $ v\cdot x$ is even for all
$x$ in the both eigenlattices of $H^n(X_+)$, and
$d(c_+)-d(c_-)=-1$ otherwise.

In the case of cubic fourfolds, we have shown the converse: the
relation $c_-\circ c_+=R_v$ turns out to be a sufficient criterion
for (coarse) adjacency of the deformation components (our original
argument in \cite{FK1} was greatly simplified in \cite{FK2}, after
Theorem~\ref{C-surj} appeared and Theorem~\ref{rsurjectivity} was
derived from it). This can be summarized as follows.

\lemma\label{wall-crossing} Real cubic fourfolds $X_+$ and $X_-$
represent adjacent coarse deformation components if and only if
there exists a lattice isomorphism $\M(X_+)\to\M(X_-)$ which
preserves the polarization classes $h_{X_\pm}\in\M(X_\pm)$ and
identifies the involution $c_-$ with the composition $R_v\circ
c_+=c_+\circ R_v$ for some $2$-root $v\in\M_\pm(X_+)$.
 Moreover,
the move leading from the deformation component of $X_+$ to that
of $X_-$ is \roster \item $R$-move, if $v\in\M_+(X_+)$ and $v\cdot
x=0\mod 2$ for all $x\in\M_+(X_+)$, \item $L$-move, if
$v\in\M_-(X_+)$, and $v\cdot x=0\mod 2$ for all $x\in\M_-(X_+)$,
 \item
$L^{-1}$-move if $v\in\M_+(X_+)$ and $v\cdot x=1\mod 2$ for some
$x\in\M_+(X_+)$,
 \item
$R^{-1}$-move if $v\in\M_-(X_+)$ and $v\cdot x=1\mod 2$ for some
 $x\in\M_-(X_+)$.
 \endroster
\qed
\endlemma

Now let us apply Theorem \ref{rsurjectivity} to the cuspidal
cubics similarly.

As it follows from analysis of the local model performed in
\ref{cusp-perturb}, there exist two ways of perturbation of a real
cuspidal hypersurface $X_0$ (of any dimension $n$): one
perturbation does not change topologically its real locus, i.e.,
yields $X_+$ with $X_+(\R)$ homeomorphic to $X_0(\R)$, and the
other perturbation adds a handle (of an appropriate index) to
$X_0(\R)\cong X_+(\R)$, and in particular, yields $X_-$ with
$b_*(X_-(\R))=b_*(X_+(\R))+2$ (or equivalently,
$d(X_-)=d(X_+)-1$).
 We explained also that
the latter perturbation shifts $X_0$ from a cuspidal stratum of
$\D(\R)\subset P_{n,d}(\R)$ inside the thin region of the cuspidal
slice of $\D(\R)$ shown on Figure~3.

Analysis of the same local model $\sum_{i=1}^p x_i^2-\sum_{j=1}^q
y_j^2+z^3+ bz +c$ of a cuspidal perturbation shows that for even
$n=p+q$, the involution $c_-\:H^n(X_-)\to H^n(X_-)$ acts in the
Milnor lattice $A_2$ of the cusp as $(-1)^q$.
 In the other words, $A_2$ is
embedded into the corresponding eigenlattice. (Note that the
involution $c_+$ in $H^n(X_+)$ interchanges some pair of roots,
$v_1,v_2$, generating $A_2$.) As it follows from Theorem
\ref{rsurjectivity}, such $A_2$-sublattice should not contain
$6$-roots (since otherwise the corresponding periods, namely the
intersection of mirrors $H_\delta$ over generating 2-roots of
$A_2$, would lie in $\Cal H_\infty$).

 In the case of cubic fourfolds, existence of such
appropriate $A_2$-sublattice yields an arithmetical criterion
allowing to detect cuspidal strata on the boundary of a
deformation component of real cubic fourfolds.

\lemma\label{cusp-lemma} Consider the wall $\Cal W$ separating a
pair of adjacent components represented by cubics $X_\pm$.
 Assume that $d(X_+)>d(X_-)$. Then $\Cal W$ contains a cuspidal
 stratum if and only if
\roster\item there exist $2$-roots $v_1,v_2\in\M_-(X_-)$ such that
 $$v_1\cdot v_2=-1 \text{ and } (v_1-v_2)\cdot
\M_-(X_-)\ne 0\mod 3,\ \text{ if \ $\Cal W$ is $R$-wall;}$$
 \item there exist
$2$-roots $v_1,v_2\in\M_+^0(X_-)$ such that $$v_1\cdot v_2=-1
\text{ and } (v_1-v_2)\cdot \M_+^0(X_-)\ne 0\mod 3,\
 \text{ if \ $\Cal W$ is $L$-wall.}$$
\endroster
\endlemma

\demo{Proof} Any $A_2$-lattice contains precisely six elements of
square six, namely, $\pm(v_1-v_2)$, $\pm(2v_1+v_2)$, and
$\pm(v_1+2v_2)$, where $v_1,v_2$ are $2$-roots forming a basis.
The two latter pairs of vectors are congruent modulo $3$ to the
first pair (since, $(2v_1+v_2)+(v_1-v_2)=3v_1$, etc.). Therefore,
under the hypothesis that $(v_1-v_2)\cdot \M^0(X_-)\ne 0\mod 3$, a
sublattice $A_2$ does not contain $6$-roots and so, by
Theorem~\ref{rsurjectivity} a generic point in $\Cal D\cap\{w\,|\,
w\cdot v_1= w\cdot v_2=0\}$ (and, as a consequence, a generic real
point in this intersection) does not belong to $\Cal H_\infty$.
Since $A_2$ is a sublattice of $\M_-(X_-)$ (of $\M_+^0(X_-)$) in
the case of an $R$-move (respectively, $L$-move), the condition
$(v_1-v_2)\cdot \M^0(X_-)\ne 0\mod 3$ can be equally stated as is
done in the Lemma.
 \qed\enddemo

\corollary\label{R-cuspidal-strata} Each of the $R$-walls, except
the ones which lead to
 $C^{10,1}$ or $C^{2,1}_I$, has a cuspidal
stratum.
\endcorollary

\proof To apply Lemma \ref{cusp-lemma} we need to embed suitably
$A_2$ to $\M_-( X_-)$. According to Tables 2 and 3, in each of the
corresponding cases the lattice $\M_-(X_-)$ contains as a direct
summand either $\langle 2\rangle\oplus U$, or $D_4$, or $E_7$, or
$E_8$. In $D_4, E_7$, and $E_8$ we may take any standard embedding
of $A_2$. In $\langle 2\rangle\oplus U$ we pick $v_1=e -u_1$ and $
v_2=u_1+u_2$, where $e$ is a generator of $\langle 2\rangle$ and
$u_1,u_2$ are standard generators of $U$.
\endproof

\section{Ramified connected sums}\label{ramified}

\subsection{Double ramified coverings and their Morse
modifications} Assume that $U$ is a compact $n$-manifold and $L$ a
codimension two submanifold coming transversely to $\partial U$
along its boundary $\partial L\subset\partial U$. Recall that
double coverings $\pi\:\til U\to U$ ramified along $L$ are
classified up to isomorphism by the characteristic class $w_1\in
H^1(U\smallsetminus L;\Z/2)$ of the restriction of $\pi$ over
$U\smallsetminus L$ (that is of the unramified part of $\pi$).
Given an element $w\in H^1(U\smallsetminus L;\Z/2)$, there exist a
covering ramified along $L$ with the characteristic class $w_1=w$
if and only if the coboundary map $H^1(U\smallsetminus L;\Z/2)\to
H^2(U,U\smallsetminus L;\Z/2)$ sends $w$ to the Thom class of the
normal bundle of $L$. As is known, the latter condition is
equivalent to a possibility to realize the class dual to $w_1$ in
$H_{n-1}(U,\partial U\cup L;\Z/2)$ by a codimension one compact
submanifold $F\subset U$ whose boundary splits as $\partial
F=L\cup F^\partial$, where $F^\partial=\partial F\cap\partial U$,
so that $\partial L=\partial F^\partial$ and $F$ is transversal to
$\partial U$ along $F^\partial$. In constructions we allow $F$ and
$F\cap \partial U$ to have real algebraic singularities outside
$L$ (then still by duality, $F$ defines a class $w_1$ as above; in
fact, Morse singularities are sufficient for our purposes) and
call such $F$ a {\it characteristic hypersurface} of the ramified
covering $\pi$ or a {\it characteristic membrane bounding} $L$.
For simplicity, {\fib we suppose from now on that $L$ and
$F^\partial$ are disjoint parts of $\partial F$.}

It may be worth recalling that given a non-singular characteristic
hypersurface $F$, a representative, $\pi\:\til U^F\to U$, of the
associated class of ramified covering is obtained from a disjoint
union $U\+U$ of two copies of $U$ by cutting them along $F$ and
then gluing together along $F$ in a cross-like fashion, so that
each side of $F$ in one copy of $U$ is identified with the
opposite side of $F$ in the other copy. Note that $F$ does not
need to be orientable, and it is sufficient that ``the sides'' are
just locally defined.

We need to treat the coverings in one-parameter families. More
precisely, to compare coverings of $U$ with coverings of
$U_\bullet=U\times[0,1]$, and {\it vice versa}. To include a
covering of $U$ into a covering of $U_\bullet$ it is necessary and
sufficient to extend $L$ to a proper codimension two submanifold
$L_\bullet$ of $U_\bullet$ and $w_1$ to a proper class of its
complement. It is for this task that characteristic hypersurfaces
are especially convenient: given an extension $L_\bullet$ of $L$
it is sufficient to extend $F$ to a characteristic membrane
$F_\bullet$ bounding $L_\bullet$.

In what follows we are concerned exclusively with a particular
situation where the projection to $[0,1]$ restricted to
$L_\bullet$ is a Morse function, so that we may consider
$L_\bullet$ as a special family $L_t$, $L_\bullet=\cup_{t\in
[0,1]}L_t\times t$, such that $L_t$ is a smooth isotopy which
experiences ambient Morse modifications at finitely many points.
As to $F_\bullet$, we restrict ourselves to continuous families,
$F_\bullet=\cup_t F_t$, where $F_t$ is a characteristic membrane
of $L_t$ for each noncritical value of $t$. (Recall that
smoothness of $F_\bullet$ is not required, and moreover, it can be
even not a topological submanifold: it is enough that it gives a
homology cycle.)

It is transparent that such a construction leads to
diffeomorphisms $\til U^{F_0}\to \til U^{F_t}$, if $L_t$ is an
isotopy. It is also a well-known fact in the knot theory, that a
Morse modification of index $q$ performed on $L_t$ yields a Morse
modification of index $q+1$ on $U^{F_t}$.
 Still, we sketch the proof below, since we could not find a reference
 suitable for our needs.

\lemma\label{r-modification} {\rm (1)} If $L_t$ is an isotopy,
then the family $F_t$ yields a continuous family of
diffeomorphisms $\til{\phi}_t\:\til U^{F_0}\to \til U^{F_t}$. If
$F_t^\partial=\oo$ for all $t\in[0,1]$, then such a family
$\til\phi_t$ is identical on $\partial\til U^{F_t}=\partial
U\+\partial U$.

{\rm (2)} Assume that there is only one critical value
$t\in(0,1)$, so that $L_1$ is obtained from $L_0$ by a Morse
modification of index $q$. Then $\til U^{F_1}$ is obtained from
$\til U^{F_0}$ by a Morse modification of index $q+1$.
\endlemma

\demo{Proof} An isotopy of the branching locus $L_t$ in (1) can be
extended to an ambient isotopy $\phi_t\:U\to U$. The latter lifts
to a family of diffeomorphisms $\til\phi_t:\til U^{F_0}\to \til
U^{F_t}$, since the characteristic classes $w_1$ of all the
coverings do match. If $F_t^\partial=\oo$, then we may choose
$\phi_t$ identical on $\partial U$ and its lifting $\til\phi_t$
identical on $\partial U\+\partial U$.

To prove (2), it is sufficient to analyze a local model of an
elementary cobordism of index $q$. Namely, it is enough to
consider the case of $U=D^n$ ($n$-ball), where $L\subset
U\times[0,1]$ is $(n-1)$-ball, and $L_t$ represent a Morse
modification of index $q$ inside $U\times[0,1]$.
 Such elementary cobordism connects standard framed spheres,
$L_0=D^{p}\times S^{q-1}\subset U\times0$ and $L_1=S^{p-1}\times
D^{q}$, $p+q=n-1$, on the bottom and the top of the cylinder
$U\times[0,1]$, and is restricted to $\partial U\times[0,1]$ as a
product-cobordism, $\partial L_t=S^{p-1}\times
S^{q-1}\subset\partial U\times t$.
 Passing to the double covers over $U\times t$ ramified along
$L_t$, we observe an index $q+1$ model cobordism between
$D^p\times S^q\times D^1$ and $S^{p}\times D^{q}\times D^1$.
\qed\enddemo

\subsection{Ramified connected sums}\label{rsum-subsection}
Assume that $F^\partial=\oo$ and $U$ is embedded into each of two
closed $n$-manifolds $X_1$ and $X_2$. Then we can remove $U$ from
each of $X_i$ and glue instead $\til U^F$ via
 the identity map of the
boundary $\partial\til U^F =\partial U\+\partial
U\subset(X_1\smallsetminus U)\cup(X_2\smallsetminus U)$. The new
manifold will be denoted $X_1\#^F X_2$ and called the {\it
ramified connected sum of $X_1$ and $X_2$}. It is easy to see that
$X_1\#^F X_2$ depends only on the embedding of $F$ into $X_1$ and
$X_2$, but not on the particular choice of its neighborhoods $U$.
Moreover, if we consider a continuous family $F_t$ like in Lemma
\ref{r-modification}, then we can apply it to obtain the
following.

\corollary\label{rsum-modification} {\rm (1)} If $L_t=\partial
F_t$ is an isotopy, then there is a continuous family of
diffeomorphisms $\psi_t\:X_1\#^{F_0} X_2\to X_1\#^{F_t} X_2$.

{\rm (2)} Assume that there is only one critical value
$t\in(0,1)$, so that $L_1$ is obtained from $L_0$ by a Morse
modification of index $q$. Then $X_1\#^{F_1} X_2$ is obtained from
$X_1\#^{F_0} X_2$ by a Morse modification of index $q+1$.\qed
\endcorollary

In what follows we will need the following particular example.
Assume that $F'=F\cup T$, where $T$ is a ``solid torus''
$D^p\times S^{q-1}$ embedded into some $n$-ball $D^n\subset
U\smallsetminus F$, $p+q=n$.

\lemma\label{add-torus} If $T$ is unknotted in $D^n$, then
$X_1\#^{F'} X_2$ is diffeomorphic to a connected sum of $X_1\#^F
X_2$ with  $(S^1\times S^{n-1})\#(S^{p}\times S^{q})$.
\endlemma

\demo{Proof} The double covering over $S^n$ ramified along an
unknotted $S^{p-1}\times S^{q-1}$ is well-known to be $S^p\times
S^q$. Note that $X_1\#^{F'} X_2$ is obtained from $X_1\#^F X_2$
after removing a pair of balls (pull-back of $D^n$) and gluing
instead the double covering over $D^n$ ramified along $\partial
T$, which is diffeomorphic to $S^p\times S^q$ with a pair of balls
removed.
 The one-handle involved is orientable (since
ramified coverings preserve orientability), and thus we just take
a connected sum of $X_1\#^F X_2$ with $(S^1\times
S^{n-1})\#(S^{p}\times S^{q})$. \qed\enddemo

\subsection{Perturbation of the union of a quadric with a
hyperplane}\label{perturb} Here we describe topologically the
result of perturbation of the union $\PP\cup \QQ\subset\Rp{n+1}$
of a real hyperplane $\PP$ and a real quadric $\QQ$.

Let $[x_0:\dots x_n:y]$ denote homogeneous coordinates in
$\Rp{n+1}$ and $\PP=\Rp{n}$ be defined by $y=0$. We assume that
$\QQ$ is defined by equation $f_2-\e y^2=0$ where
$f_2=f_2(x_0,\dots,x_n)$ is a non-degenerate quadratic form and
$\e>0$ is a fixed parameter.
 The intersection $\VV=\PP\cap \QQ=\{f_2(x)=0\}$ gives a splitting
$\PP=P_+\cup P_-$, where $P_\pm=\{[x]\in \PP\,|\,\pm
f_2(x)\ge0\}$. If $f_2$ has signature $(p,q)$, $p+q=n+1$, then
$P_+$ is a tubular neighborhood of $\Rp{p-1}$ and $P_-$ the
tubular neighborhood of a complementary $\Rp{q-1}$. Forgetting the
coordinate $y$ gives the projection $p\:\QQ\to P_+$ which is
obtained from the orientation double covering $\til p\:\til P_+\to
P_+$ by identifying the pairs of points $\til p^{-1}(x)$ along the
boundary.

\lemma\label{odd-degree-locus} The locus $F_+=P_+\cap\{f=0\}$ of
any non-singular polynomial $f$ of odd degree is a characteristic
hypersurface of the double covering $\til p\:\til P_+\to P_+$.
\endlemma

\demo{Proof} An odd degree hypersurface $\{f=0\}$ in $\PP$ is dual
to the generator of $H^1(\PP;\Z/2)$ and so is characteristic for
the orientation covering $\til \PP\to \PP$. Restricting the latter
over $P_+$ we conclude that $F_+$ is characteristic for $\til
P_+\to P_+$. \qed\enddemo

By a perturbation of $\PP\cup \QQ$ via $f_3=f_3(x_0,\dots,x_n,y)$
we will mean a real cubic $\ZZ$ defined by $y(f_2-\e y^2)+\delta
f_3=0$, where $0<\delta<\!<\e$. Such a cubic is non-singular
provided that $P$, $Q$ and $f_3=0$ (or equivalently, $f_2(x)=0$
and $f_3(x,0)=0$ inside $P$) intersect transversally.

Denote this intersection by $L$ and assume the above
transversality. Then, the intersection of $\{f_3=0\}$ with $P_+$
is a characteristic membrane bounding $L$.

Consider the double $D(P_+)$ of $P_+$, which is obtained from two
copies of $P_+$ by gluing them together along the boundary $\VV$.
Let $U$ denote a neighborhood of $P_+$ obtained by adding to it a
collar $\VV\times[0,1)$ of its boundary, so that $F\subset
P_+\subset\Int(U)$. Note that $P_+$ lies both in $P$ and in
$D(P_+)$, and we may extend these embeddings of $P_+$ to
embeddings $U\subset P$ and $U\subset D(P_+)$.

Our next target is the following proposition, which will play a
key role in the forthcoming analysis of cubic fourfolds.

\proposition\label{rsum} If the quadratic form $f_2(x)$ is
non-degenerate and the intersection of $f_2(x)=0$ and $f_3(x,0)=0$
in $P$ is transversal, then for any $0<\delta<\!<\e$ the real
locus $Z(\R)$ of the cubic $Z=\{y(f_2-\e y^2)+\delta f_3=0\}$ is
diffeomorphic to $\PP\#^F D(P_+)$.
\endproposition

\subsection{The perturbation neck}\label{subsection-neck}
We start proving Proposition \ref{rsum} with checking certain
properties of the perturbations of $\PP\cup \QQ\subset\Rp{n+1}$ in
a somewhat more general setting, assuming that $\PP$ and $\QQ$ are
transverse real non-singular hypersurfaces defined by polynomials
$f$ and $g$ of degrees $\ell$ and $m$. As before, we let
$\VV=\PP\cap \QQ$ and denote by $\ZZ\subset\Rp{n+1}$ the result of
a perturbation described by a polynomial $fg+\delta h$,
$0<\delta<\!<1$, assuming in addition that $P,Q,$ and $h=0$
intersect transversally. Thus, their intersection, which we denote
$\LL$, is non-singular.

It is well known (and trivial) that topologically a small
perturbation of $\PP\cup \QQ$ is localized in a tubular
neighborhood $N\subset\Rp{n+1}$ of $\VV$. Namely, the tubular
neighborhoods $N_P=N\cap \PP$ and $N_Q=N\cap \QQ$ of $\VV$ are
removed from $\PP$ and $\QQ$ and are replaced by a certain {\it
perturbation neck} $N_{P,Q}\subset N$, which has also a projection
$p_{P,Q}\:N_{P,Q}\to\VV$.
 Namely, the fibers  $I_v^P$ and $I_v^Q$ of the projections $N_P\to \VV$ and
 $N_Q\to \VV$ over $v\in \VV$ look topologically like line segments
 intersecting in the middle points.
The cross-like fiber $I_v^P\cup I_v^Q$ of $(\PP\cup \QQ)\to \VV$
is perturbed into a fiber of $p_{P,Q}$, which looks generically as
a pair of arcs $I\+I$ connecting the two endpoints $\partial
I_s^P$ with the two endpoints $\partial I_s^Q$, and the two
possible ways of such connection alternate as we cross the locus
$\LL=\VV\cap\{h=0\}$. Over $\LL$ the fibers of $p_{P,Q}$ are
non-generic, namely, they remain $I_v^P\cup I_v^Q$.

\midinsert \topcaption{Figure 4. Factorization of a perturbation
neck}
\endcaption\epsfbox{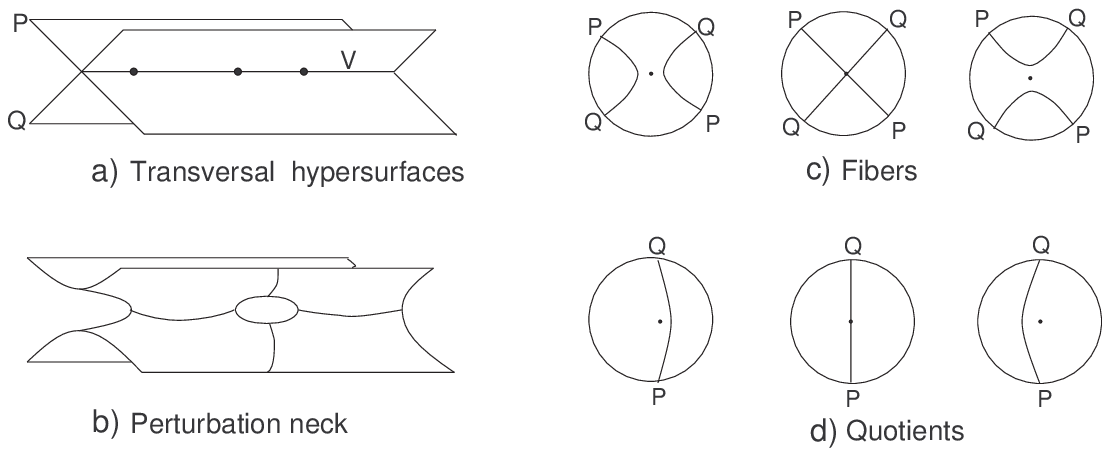}
\botcaption{\eightpoint \noindent a) Hypersurfaces $P\cup Q$,
where $L\subset V=P\cap Q$ is marked by dots. b) The corresponding
perturbation neck
 $N_{P,Q}$.
c) Fibers of $p_{P,Q}\:N_{P,Q}\to V$. The central singular fiber
$I_v^P\cup I_v^Q$ is over $v\in L$; it is perturbed in two ways as
$v$ shifts from $L$. d) Passing to the quotient by $\tau$ (central
symmetry) makes all the fibers homeomorphic to $I$.}\endcaption
\endinsert

Next we observe existence of an involution $\tau\:N_{P,Q}\to
N_{P,Q}$ which preserves the fibers and act on each fiber as the
central symmetry. The quotient space $N_{P,Q}/\tau$ is fibred over
$\VV$ trivially with a fiber $I$, and so can be identified with
$\VV\times[-1,1]$. This implies that the quotient map
$q:N_{P,Q}\to \VV\times[-1,1]$ has to be a ramified covering.

\lemma\label{neck} The quotient map $q\:N_{P,Q}\to
\VV\times[-1,1]$ is a double covering ramified along $\LL\times0$.
\qed\endlemma

Let us choose the identification $N_{P,Q}/\tau=\VV\times[-1,1]$ so
that $\partial N_P/\tau$ and $\partial N_Q/\tau$ are identified
with $\VV\times-1$ and $\VV\times1$ respectively.

The lemma below describes the characteristic hypersurface $F$ of
the double covering $q$ under an additional assumption that one of
the degrees, $m=\deg(Q)$, is even. This assumption implies that
$\VV$ is co-orientable in $\PP$ and thus its tubular neighborhood
$N_P$ can be identified with $N_{P,Q}/\tau=\VV\times[-1,1]$.

\lemma\label{r-neck} Assume that $m$ is even. Then the
intersection $F=(\VV\times[0,1])\cap\{h=0\}$ is the characteristic
hypersurface of the ramified covering $q$, with respect to the
identification of the quotient space $N_{P,Q}/\tau$ with
$N_P=\VV\times[-1,1]$, as described above.
\endlemma

\demo{Proof} Note that $q$ has the ramification locus $\LL=F\cap
(\VV\times0)=\{f=g=h=0\}$, as required by Lemma \ref{neck}. The
restriction of $q$ over $\VV\times{-1}$ is trivial which
corresponds to the characteristic cycle
$\oo=F\cap(\VV\times{-1})$. This implies that
$\VV\times[0,1]\cap\{h=0\}$ is the characteristic cycle of $q$
because $H^1(\VV\times[-1,1],\VV\times{-1};\Z/2)=0$. \qed\enddemo

Now, fixing $f$ and $g$ let us assume that $h$ varies continuously
so that $\LL$ experiences a single Morse modification.

\corollary\label{general-deform-perturb} If $h$ varies so that
$\LL$ experiences a Morse modification of index $q$, then $\ZZ$
experiences a Morse modification of index $q+1$.
\endcorollary

\demo{Proof} The pieces $\PP\smallsetminus N_P$ and
$\QQ\smallsetminus N_Q$ do not change in process of variation
 of $\ZZ$, and the
perturbation neck $N_{P,Q}$ experiences a Morse modification of
index $q+1$, as follows from Lemmas \ref{r-neck} and
\ref{r-modification}. \qed\enddemo

\subsection{Proof of Proposition \ref{rsum}}
According to \ref{subsection-neck}, a perturbation of $\PP\cup
\QQ$ gives $\ZZ=(\PP\smallsetminus N_P)\cup(\QQ\smallsetminus
N_Q)\cup N_{P,Q}$, whereas $\PP\#^FD(P_+)=(\PP\smallsetminus
U)\cup(D(P_+)\smallsetminus U)\cup\til U^F$.
 Note that $\PP\smallsetminus N_P$ splits into a pair of
components, namely, $\PP\smallsetminus U$ and
$D(P_+)\smallsetminus U$, which are diffeomorphic to $P_-$ and
$P_+$ respectively. So, the proof of Proposition \ref{rsum} is
reduced to the following statement.

 \lemma The part $(\QQ\smallsetminus N_Q)\cup N_{P,Q}\subset \ZZ$
is diffeomorphic to $\til U^F$.
 \endlemma

\demo{Proof} Recall that $U$ is $P_+$ enlarged with a collar in
$P_-$, so we can take $U=(P_+\smallsetminus N_P)\cup N_P$.

Then $\til U^F$ splits into a union of two pieces, one of which is
the pull-back of $N_P$ and the pull-back of $\Cl(P_+\smallsetminus
N_P)$. Lemma \ref{r-neck} says that the first piece is just
$N_{P,Q}$. The other piece is obviously diffeomorphic to the
pull-back of $P_+$. It can be identified with $\til P_+$, because
of Lemma \ref{odd-degree-locus}.
 To complete the proof it is left to recall
that $\QQ\smallsetminus N_Q$ is diffeomorphic to $\til
P_+$.\qed\enddemo

\subsection{Perturbation of nodal cubics}\label{perturb-nodal}
Like in \ref{perturb}, let us consider a hyperplane $\PP=\{y=0\}$
and a quadric $\QQ=\{f_2(x)=\e y^2\}$ in $\Rp{n+1}$ with
coordinates $x=(x_0,\dots,x_n)$ and $y$, where $f_2$ is a real
non-degenerate quadratic form. But now let us choose a cubic form
$f_3=f_3(x)$ independent of $y$. Suppose like before that the
intersection $Y=\{f_2=f_3=0\}$ is transversal (then $Y$ is a
non-singular $(n-1)$-fold), and denote by $X_{\e,\delta}$ the
cubic $n$-fold $\{y(f_2-\e y^2)+\delta f_3=0\}$. In \ref{perturb},
under the assumption that $0<\delta<\!<\e$, we treated
 $X_{\e,\delta}(\R)$ as a perturbation of
$\PP\cup\QQ$. One can also view $X_{\e,\delta}$ as a perturbation
of a nodal cubic $\{yf_2+\delta f_3=0\}$ supposing that
$0<\e<\!<\delta$. The following Lemma shows, in particular, that
both perturbations yield the same deformation type.

\lemma\label{nonsingularity} For any $f_2(x)$, $f_3(x)$ like
above, there exists $\varkappa>0$ such that for any choice of
positive $\e,\delta<\varkappa$ the locus
$X_{\e,\delta}=\{y(f_2(x)-\e y^2)+\delta f_3(x)=0\}$ in $\Rp{n+1}$
is non-singular.
\endlemma

\proof Under the quasi-homogeneous change of variables
$(x,y)\mapsto (tx,t^{-2} y)$ the polynomial $y(f_2(x)-\e
y^2)+\delta f_3(x)$ is transformed into
$$t^{-2} y(f_2(tx)-\e (t^{-2} y)^2)+\delta f_3(tx)=y(f_2(x)-\e t^{-4}
y^2)+\delta t^3 f_3(x),$$ so that the plane of parameters $\e,
\delta$ becomes foliated by curves $\e^3\delta^4=$\rm{const} each
representing projectively equivalent $X_{\e,\delta}$. Hence, the
statement follows from the non-singularity of $X_{1,\delta }$ for
small $\delta >0$.
\endproof

From now on we stick to the case of our interest, $n=4$.

\proposition\label{k3-k4} Consider a vertex $v_{K3}$ of the
K3-graph and a vertex $v_{K4}$ of the K4-graph which correspond to
each other. Assume that a non-degenerate quadric $\{f_2=0\}$ and a
cubic $\{f_3=0\}$ intersect transversally, their intersection
$\{f_2=f_3=0\}$ represents the vertex $v_{K3}$, and the sign of
$f_2$ is chosen so that $\chi(P_+)=1$, where $P_+=\{x\in
\PP\,|\,f_2(X)\ge0\}$. Then there exists $\varkappa>0$ such that
for any choice of positive $\e,\delta<\varkappa$ the cubic
fourfold $X_{\e,\delta}=\{y(f_2-\e y^2)+\delta f_3=0\}$ is
non-singular and represent vertex $v_{K4}$.
\endproposition

\demo{Proof} The nodal cubic fourfolds $X_{0,\delta}$ defined by
$yf_2+\delta f_3$ are in central projection correspondence with
$Y$  and represent an edge in the K4-graph. By theorem
\ref{k3-k4correspondence}, it is the leftmost endpoint (providing
a smaller value of the coordinate $r$) of the edge which
corresponds to  the vertex $v_{K3}$ representing $Y$. The two
endpoints are given by two perturbations, $\{y(f_2-\e y^2)+\delta
f_3\}$ and $\{y(f_2+\e y^2)+\delta f_3\}$, of $X_{0,\delta}$. The
relation $2r=b_*(X)-\chi(X(\R))-4$ (see \ref{K4}) implies that the
vertex $v_{K4}$ is represented by that perturbation which yields a
cubic fourfold with a greater value of $\chi(X(\R))$. As it
follows from Proposition \ref{rsum},
$\chi(X_{\e,\delta}(\R))=\chi(\PP)+2\chi(P_+)-\chi(\LL)$, which is
greater than the Euler characteristic for another perturbation,
$\chi(\PP)+2\chi(P_-)-\chi(\LL)$, since by our assumption
$\chi(P_+)=1$ and thus $\chi(P_-)=0$.

Finally, we apply Lemma \ref{nonsingularity}.\qed\enddemo

\corollary\label{real-locus-k3-k4} Under the assumption of
Proposition \ref{k3-k4}, the real locus of a cubic fourfold of the
type $v_{K4}$ is diffeomorphic to $\PP\#^F D(P_+)$, where
$F=\{x\in P_+\,|\,f_3(x)=0\}$.
\endcorollary

\demo{Proof} Choosing $\delta<\!<\e$ we can view $Y_{\e,\delta}$
as a perturbation of $\{y(f_2-\e y^2)\}$ and apply Proposition
\ref{rsum}. \qed\enddemo

Let us fix $f_2$ and continuously vary $f_3$, so that the real
K3-surfaces $Y_{t}$, defined by $f_{3,t}=f_2=0$, $t\in[-1,1]$, are
non-singular for $t\ne0$, and $Y_{0}$ is nodal, so that $Y_{\pm
t}$, where $t>0$, represent a pair of adjacent vertices
$v^{\pm}_{K3}$ of the K3-graph. Let $X_{\pm}$ be cubic fourfolds
representing the corresponding vertices $v^{\pm}_{K4}$ of the
K4-graph.

\proposition\label{index+1} If $Y_{t}(\R), t<0,$ experiences  at
$t=0$ a Morse modification of index $q$, then $X_{+}(\R)$ is
obtained from $X_{-}(\R)$ by a Morse modification of index $q+1$.
\endproposition

\demo{Proof} According to Corollary \ref{real-locus-k3-k4}, we
have $X_{\pm}(\R)=\PP\#^{F_{\pm}} D(P_+)$, where
$F_{\pm}=P_+\cap\{f_{3,\pm 1}=0\}$. Now it remains to apply
Corollary \ref{rsum-modification}.
 \qed\enddemo

\section{Proof of the main theorem}\label{main-thm}

\subsection{The proof in the general case}
We have already proved Theorem \ref{main-theorem} in the cases
$C^{1,0}_I$ and $C^{0,0}$, see \ref{immediate}. As the next step,
let us consider the real cubic fourfolds, $X$, representing the
vertices of $\Gamma_{K4}$ which are obtained from $C^{0,0}$ by
L-moves (i.e., lying on the upper-left side of  $\Gamma_{K4}$).

\lemma\label{left-top} If a cubic $X$ represents $C^{i,0}$, $i\ge
0$, or $C^{i,0}_I$, $i> 1$, then $X(\R)$ is diffeomorphic to
$\Rp4\#i(S^2\times S^2)$.
\endlemma

\demo{Proof} We prove this lemma by induction in $i$. The base
case, $i=0$, is already established in Corollary \ref{basecase}.
Consider a cubic fourfold $X'$ representing $C^{i-1,0}$, $i>0$,
and assume that $X$ represents either $C^{i,0}$ or $C^{i,0}_I\ne
C^{1,0}_I$. By Lemma \ref{LRindex}(3), $X(\R)$ is obtained from
$X'(\R)$ by attaching a 2-handle, and by Lemma \ref{LRparity}(2)
the core circles of this handle is null-homologous. By inductive
assumption, $X'(\R)= \Rp4\#(i-1)(S^2\times S^2)$ has $\pi_1=\Z/2$,
and thus the core circle is contractible. Attaching of a
$2$-handle along a contractible circle is equivalent to taking a
connected sum either with $S^2\times S^2$, or with
$\Cp2\#\overline{\Cp2}$, depending on the framing of the
$2$-handle.
 The second option is
impossible because all cubic fourfolds have $w_2(X(\R))=0$, and
thus, the orientation covering space of $X_\R$ must have even
intersection form. \qed\enddemo

Our next aim is to apply $R$-moves to the cubics treated above.

\lemma\label{Rmove} If the R-wall separating the cubics $X'$ of
type $C^{i,j-1}$, $j\ge1$, from the cubics $X$ of type $C^{i,j}$
or $C^{i,j}_I$ contains a cuspidal stratum, then $X_\R$ is
diffeomorphic to $X_\R'\#(S^1\times S^3)$.
\endlemma

\demo{Proof} By corollary \ref{cusp-stratum}(1), one of the
facet-strata adjacent to the cuspidal stratum has index $1$, and
so $X(\R)$ is obtained by adding a $1$-handle to $X'(\R)$. The
latter is connected by Lemma \ref{disconnected} and Corollary
\ref{basecase}. It is also non-orientable, so, adding $1$-handle
 means taking a connected sum with $S^1\times S^3$.
\qed\enddemo

This allows us to deduce the main theorem \ref{main-theorem} in
all but a few cases.

\corollary\label{general-case} If a cubic $X$ belongs to the type
$C^{i,j}$ or $C^{i,j}_I$ different from $C^{1,0}_I$, $C^{10,1}$,
and $C^{2,1}_I$, then $X_\R$ is diffeomorphic to
$\Rp4\#i(S^2\times S^2)\#j(S^1\times S^3)$.
\endcorollary

\demo{Proof} Lemma \ref{left-top} covers the cases with $j=0$.
Vertices $C^{i,j}$ and $C^{i,j}_I$, $j>0$ can be reached from
$C^{i,0}$ by R-moves. According to Corollary
\ref{R-cuspidal-strata}, the assumption of Lemma \ref{Rmove} is
satisfied unless $X$ belongs to the type $C^{10,1}$, or
$C^{2,1}_I$. So, applying Lemma \ref{Rmove} $j$ times we obtain
the given description of $X(\R)$. (Note that the exceptional
vertices are terminal and thus are not obstacles in the sequence
of R-moves.)
 \qed\enddemo

\subsection{The case of $C^{2,1}_I$}
 Like in \ref{perturb}, we
start with a perturbation of a reducible real cubic fourfold
$\PP\cup\QQ$, where quadric $\QQ$ is defined by $f_2-\e y^2$. Now
we specify $f_2$ so that the region $P_+$ is a tubular
neighborhood of $P^2(\R)=\{[x_0\:\!x_1\:\!x_2\:\!0\:\!0]\}\subset
P$, for example, we may take $f_2=x_0^2+x_1^2+x_2^2-x_3^2-x_4^2$.
As a perturbation term $f_3$ we pick a degree 3 homogeneous
polynomial in three variables $x_0,x_1,x_2$ (thus, independent of
$x_3,x_4,y$) such that the curve defined in $P^2(\R)$ by equation
$f_3=0$ is non-singular, which insures in particular the
transversality of intersection between $f_3=0$ and $f_2=0$ in
$P^5(\R)$.

Note that $F=P_+\cap\{f_3=0\}$ is the pull-back of this curve
under the tubular neighborhood projection $q\:P_+\to P^2(\R)$
(that is the projection which forgets the coordinates $x_3,x_4$).
Therefore, under above choices, $F$ is diffeomorphic to the
product of the curve by $D^2$.

A real plane cubic curve may have one or two connected components,
and we need to consider the both possibilities. Let us denote by
$f_3^{(k)}$, $k=1,2$, a polynomial defining a real nonsingular
cubic curve with $k$ component. Then, $F^{(k)}=P_+\cap
\{f_3^{(k)}=0\}$ becomes a solid torus for $k=1$ and a pair of
solid tori for $k=2$.

Following Section \ref{perturb}, consider real cubic fourfolds
$X^{(k)} (k=1,2 )$ defined by equation $y(f_2-\e y^2)+\delta
f_3^{(k)}=0$, where we pick positive $\e$ and $\delta$ to be
smaller than the constant $\varkappa$ provided by Lemma
\ref{nonsingularity}.

\lemma\label{twotypes} $X^{(1)}$ and $X^{(2)}$ belong to the types
$C^{1,0}$ and $C^{2,1}_I$ respectively.
\endlemma

\demo{Proof} The real locus of the K3-surface $f_2=f_3^{(k)}=0$ is
$\partial F^{(k)}$, that is $k$ copies of a torus. Such a
K3-surface is known to be of type $C^{1,0}$ if $k=1$ and
$C^{2,1}_I$ if $k=2$ (see for example the survey \cite{DK}).
 Since $\chi(P_+)=1$, Proposition \ref{k3-k4}
implies that the cubic fourfold $X^{(k)}$ belongs to the
corresponding type, as is stated. \qed\enddemo

\lemma\label{add-torus-in-K3} $X^{(2)}(\R)=X^{(1)}(\R)\#(S^2\times
S^2)\#(S^1\times S^3)$.
\endlemma

\demo{Proof} One of the real components of the two-component cubic
curve $f_3^{(2)}=0$ is non contractible and isotopic to a real
line; the other component bounds a disc, $D\subset P^2(\R)$. Let
$T_\infty$,  $T_0$ denote the corresponding solid torus components
of $F^{(2)}$. The component $T_0$ is an unknotted torus contained
in the 4-disc $q^{-1}(D)$ disjoint from $T_\infty$. On the other
hand, the only real component of the curve $f_3^{(1)}=0$ is
isotopic in $P^2(\R)$ to a real line, and therefore $F^{(1)}$ is
isotopic to $T_\infty$. So, we can apply Lemma \ref{add-torus} to
compare $X^{(2)}=P\#^{F^{(2)}}D(Q^+)$ with
$X^{(1)}=P\#^{F^{(1)}}D(Q^+)$ and to conclude that
$X^{(2)}(\R)=X^{(1)}(\R)\#(S^2\times S^2)\#(S^1\times S^3)$.
\qed\enddemo

\corollary\label{firstcase} If a cubic $X$ belongs to the type
$C^{2,1}_I$, then $X(\R)$ is diffeomorphic to $\Rp4\#2(S^2\times
S^2)\#(S^1\times S^3)$.
\endcorollary

\demo{Proof} According to Lemma \ref{twotypes} the cubic $X^{(1)}$
belongs to type $C^{1,0}$. So, by Lemma \ref{left-top}, its real
part $X^{(1)}(\R)$ is diffeomorphic to $\Rp4\#(S^2\times S^2)$.
Thus, it remains to apply Lemma \ref{add-torus-in-K3}.
\qed\enddemo

\subsection{The case of $C^{10,1}$}

\lemma\label{K3-family} There exist a real non-degenerate
homogeneous quadratic polynomial $f_2=f_2(x_0,\dots,x_4)$ and a
continuous family of real cubic homogeneous polynomials
$f_{3,t}=f_{3,t}(x_0,\dots,x_4)$, $t\in[-1,1]$, such that the
complete intersections $Y_t\subset P^4$, $f_2=f_{3,t}=0$, are
K3-surfaces which are non-singular of type $C^{10,0}$ for $t<0$,
and non-singular of type $C^{10,1}$ for $t>0$, while $Y_0$ has a
nodal singularity.
\endlemma

\demo{Proof} According to the K3-graph, the $6$-polarized
K3-surfaces of type $C^{10,0}$ are adjacent to the K3-surfaces of
type $C^{10,1}$. Therefore, we can find a path which leads from
one component to another and intersects the wall between them at
its non-singular point. Let $f_{2,t}=f_{3,t}=0$, $t\in[-1,1]$, be
such a path, so that $t=0$ corresponds to a nodal K3-surface. For
a generic choice of such a path
 the quadric $f_{2,0}$ is non-singular and, thus, the quadrics
 $f_{2,t}$ are
also non-singular for sufficiently small $|t|>0$.
 Recall that deformation equivalent real quadrics are actually
projective equivalent, so, applying suitable projective
transformations we can turn $f_{2,t}$ into a constant family of
quadrics, $f_2$. \qed\enddemo

The topological types of $X(\R)$ for K3-surfaces $X$ are well
known (see, for example, the survey \cite{DK}). In particular, in
the case of K3-surfaces of types $C^{10,0}$ and $C^{10,1}$, the
locus $X(\R)$ is homeomorphic respectively to $S_{10}$ and
$S_{10}\+S^2$, where $S_{10}$ stands for an orientable surface of
genus $10$.
 Note also that a Morse modification of a surface
 which brings about a new component $S^2$
can be either of index $0$ (a birth of sphere), or $2$ (splitting
off a sphere from another component). Let us show that in our case
the second possibility is not realizable.

\lemma\label{only-collaps} Assume that $Y_t$, $t\in[-1,1]$ is a
family of real K3-surfaces which represent type $C^{10,0}$ for
$t<0$, type $C^{10,1}$ for $t>0$, and experience a Morse
modification at $t=0$. Then the  Morse modification is of index
$0$.
\endlemma

\demo{Proof} The component $S^2$ of $Y_1(\R)$ represents a class
$s$ in the K3-lattice $\LLL=H_2(Y_1)$. This class belongs to the
eigen-sublattice $\LLL_+$ since $S^2$ is $c$-invariant, and it has
square $s^2=-2$ since the normal bundle to the real locus is
isomorphic to the cotangent bundle. The vanishing cycle $v\in
\LLL$ of the degeneration $Y_1\to Y_0$ also has square $v^2=-2$
and, since the Morse index is even, the vanishing cycle also
belongs to $\LLL_+$. If $Y_t(\R)$ experiences a Morse modification
of index $2$, then $s\cdot v=\pm1$, so that $s$ and $v$ span a
sublattice $A_2$ in $\LLL_+$. On the other hand, the lattice
$\LLL_+$ is isomorphic to $U$ and, hence, does not contain $A_2$.
\qed\enddemo

\corollary\label{facet1} The wall separating the deformation
components $C^{10,0}$ and $C^{10,1}$ in the space of real cubic
fourfolds contains a facet of index $1$.
\endcorollary

\demo{Proof} Consider a family of real cubic fourfolds $X_t$
obtained by perturbation of $\PP\cup\QQ$, where $P=P^4(\R)$ and
$\QQ=\{f_2-y^2=0\}$, with a perturbation term $f_{3,t}$ given by
Lemma \ref{K3-family}. Namely, the family defined by
$y(f_2-y^2)+\delta f_{3,t}$ for a suitable $0<\delta<\!<1$. By
Lemma \ref{only-collaps}, the family $Y_t(\R)$ experiences a Morse
modification of index $0$ at $t=0$. Proposition \ref{index+1}
implies that $X_t(\R)$ experiences a Morse modification of index
$1$.
 \qed\enddemo

\corollary If a cubic $X$ belongs to type $C^{10,1}$, then $X(\R)$
is diffeomorphic to $\Rp4\#10(S^2\times S^2)\#(S^1\times S^3)$.
\endcorollary

\demo{Proof} The class $C^{10,0}$ has a representative $X'$ with
$X'(\R)=\Rp4\#10(S^2\times S^2)$ as we know from Corollary
\ref{general-case}. By Corollary \ref{facet1}, our $X(\R)$ can be
obtained from $X'(\R)$ by adding $1$-handle, i.e., by taking
connected sum with $S^1\times S^3$.
 \qed\enddemo

\section{Concluding remarks}\label{remarks}

\subsection{Topology of cubic threefolds}\label{3fold-subsection}
For non-singular real cubic threefolds $X$ representing seven
deformation classes among the eight existing classes with
connected real part
 $X(\R)$, the  topological
type of $X(\R)$ was determined by V.~Krasnov \cite{Kr2}. He proved
that for all  of them the real part is diffeomorphic to a
connected sum of $\Rp3$ with some number of $S^1\times S^2$. Here,
we show that the situation is different for the remaining
deformation class, that is $\BB(1)_I'$ in Krasnov's notation.

To describe the topological type of $X(\R)$ for cubic threefolds
of class $\BB(1)_I'$, we start with a construction of $X(\R)$
following the lines of Section \ref{perturb}. Namely, we consider
the cubics $X=X_{\e,\delta}=\{y(f_2-\e y^2)+\delta f_3=0\}$
obtained by a perturbation of $\PP\cup\QQ$, where $\PP=\Rp3$ is a
hyperplane $y=0$ in $\Rp4$ with coordinates
${[x_0\:\!x_1\:\!x_2\:\!x_3\:\!y]}$, and $Q\subset\Rp4$ is a real
non-singular quadric threefold, $\{f_2-\e y^2=0\}$, such that
intersection $\VV=\PP\cap\QQ$ is  transversal.
 In addition, we
suppose that the quadric surface $\VV\subset\Rp3$ is a hyperboloid
and denote by $\ell_1$, $\ell_2$ two real lines representing the
two families of real generators of this hyperboloid.
Topologically, $\VV$ is a torus which splits $\PP$ into a pair of
solid tori $P_\pm=\{\pm f_2\ge0\}$ with the meridians $m_\pm$
homologous to $[\ell_1]\pm[\ell_2]\in H_1(\VV)$. Let us orient
$\VV$ and $\ell_i$ so that $\ell_1\circ\ell_2=1$ and consider the
curve $\LL=\{f_2=f_3=0\}$ on torus $\VV$ (as in Section
\ref{perturb}, we suppose that the quadric $f_2=0$ and the cubic
$f_3=0$ are transversal to each other). B\'ezout theorem implies
that $\LL$ is homologous to $p[\ell_1]+q[\ell_2]\in H_1(\VV)$ with
$p,q\in\{\pm1,\pm3\}$, and one can easily realize all such pairs
of $p$ and $q$ by suitable choices of $f_3$. In particular, we may
consider $[\LL]=[\ell_2]-3[\ell_1]$, represented by a
spiral-shaped single-component curve $\LL$.

\theorem\label{spiral} Let $\VV$ and $\LL$ be as above. Assume
that $L$ is homologous to $[\ell_2]-3[\ell_1]\in H_1(\VV)$. Then
$X=\{y(f_2-\e y^2)+\delta f_3=0\}$ represents the deformation
class $\BB(1)_I'$ and its real part $X(\R)$ is the Seifert
manifold whose link diagram is the rightmost diagram on Figure 6.
\endtheorem

In the Figures and in the proof, we follow traditional Kirby
calculus notation and terminology, see \cite{Kirby}.

\demo{Proof} Let $T_1$ denote an abstract solid torus bounded by
$\VV$ with the meridian $\ell_1$. Since $m_+\circ\ell_1=-1$, the
union $P_+\cup T_1$ (with the common boundary $\VV$) forms a
$3$-sphere. We orient $S^3$ so that the inherited orientation on
$P_+$ restricts to the given one on $V$.

\midinsert \topcaption{Figure 5. Ramified covering along $L$}
\endcaption\epsfbox{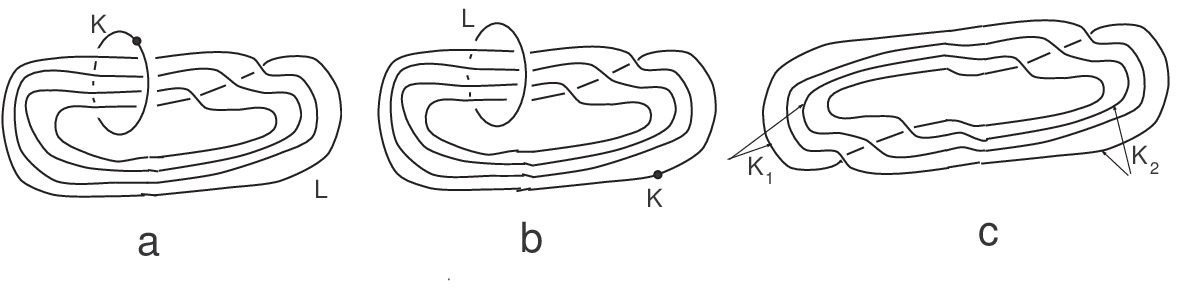}
\endinsert

On Figure 5a the core circle of $T_1$ is presented by the dotted
unknot $K$. The exterior of $K$ presents $P_+$, where $L$ looks
like a $(4,1)$-torus knot, because a positive basis of $\VV$ is
formed by classes $[m_+]$ and $-[\ell_1]$, and
$[L]=[\ell_2]-3[\ell_1]=[m_+]-4[\ell_1]$.

Note that $m_-$ can be seen as a $(2,1)$-torus knot, since
$[m_-]=-[m_+]+2[\ell_1]$. The latter gives framing $(-2)$ to $K$,
which means that $P$ is obtained from our $3$-sphere by a Dehn
$(-2)$-surgery along $K$. Note also that the double $D(P_+)$ is
obtained by $0$-surgery along $K$.

Furthermore, the Seifert surface $F=P_+\cap\{f_3=0\}$ is disjoint
from $K$, which implies that for a neighborhood $U\supset P_+$
chosen like in \ref{perturb}, our ramified double covering $\til
U^F\to U$ can be extended to the (unique) double covering
$\pi\:M_L\to S^3$ ramified along $L$, where $M_L$ it also a
$3$-sphere, since $L$ is an unknot. This implies that
$X(\R)=P\#^FD(P_+)$ is obtained from $M_L$ by a Dehn surgery along
the $2$-component link $K_1\cup K_2=\pi^{-1}(K)$ with certain
framings $n_1,n_2$ inherited from the framings $-2$ and $0$ of $K$
in $S^3$.

Note that $K$ and $L$ can be interchanged by an isotopy, and the
diagram on Figure 5b presents the same link $K\cup L$. Now $K$
looks like a $(4,1)$-torus knot, and its pre-image in $M_L=S^3$
can be seen on Figure 5c as a $(4,2)$-torus link, whose components
$K_1$, $K_2$ are $(2,1)$-torus knots. The pull-back of the torus
framing of $K$ gives the torus framing on each of $K_i$. Since the
torus framing of a $(p,q)$-torus knot is $pq$, we see that a
framing $n$ of $K$ gives framings $(n-2)$ on $K_i$, and
 thus $n_1=-4$ and $n_2=-2$. So $X(\R)$ can be described by
the leftmost link diagram on Figure 6.

\midinsert \topcaption{Figure 6. Framed link diagrams of $X(\R)$}
\endcaption\epsfbox{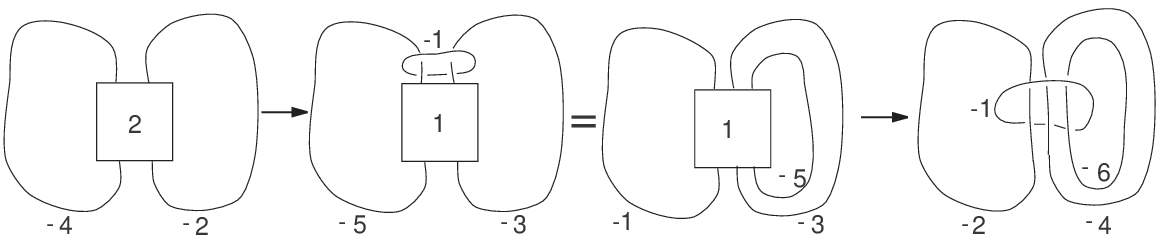}
\endinsert

The components $K_1$ and $K_2$ can be unlinked by a blowup move (a
Kirby move introducing an additional $(-1)$-components), as is
shown on Figure 6, and we obtain after another blowup move the
rightmost link diagram, which describes a Seifert 3-manifold.

The 3-manifold $X(\R)$ is obviously connected and is not
homeomorphic to $\Rp3$ with handles: it can be distinguished for
instance by its fundamental group $\pi_1(X(\R))=\la
a,b,c\,|\,a^2=b^4=c^6=abc\ra$, and even by its homology
$H_1(X(\R))=\Z/2+\Z/2$ (see for example the presentation of
$\pi_1$ given in \cite{Orlik}). So, the deformation class of our
cubic $X$ differs from the ones analyzed in \cite{Kr2}, and so has
to be the class $\BB(1)_I'$.
 \qed\enddemo

\rk{Remark} Note that Figure 6 describes the link of the
singularity $xy(x^3+y^2)+z^2=0$ known as a weighted homogeneous
singularity of type $Z_{12}$. It is one of Arnold's fourteen
exceptional singularities of modality one (its
 Dolgachev numbers $(p,q,r)$ and the Gabrielov numbers
 $(p',q',r')$ both are equal to $(2,4,6)$).
\endrk

\subsection{Cubics-handles in higher dimensions}
In higher dimensions, to construct non-singular real $n$-fold
cubics $X\subset P^{n+1}$ one may follow the same scheme as in the
proof of the main theorem: start with the case $X(\R)=\Rp{n}$ and
then make it ``develop'' and acquire more and more various handles
due to crossing of the facets of the discriminant by special moves
in the parameter space. Here, let us indicate a few first steps on
this way.

According to Lemma \ref{disconnected}, the cubics $X$ with
disconnected real part (which implies $X(\R)=\Rp{n}\+S^n$) form a
single deformation component, furthermore, this deformation
component, which we denote below by $C_0$, has only one adjacent
component, which is formed by cubics $X$ with $X(\R)=\Rp{n}$. Let
us call the latter component {\it root component} and denote it
$\Cr$. A special role of the root component is due to the
following property.

\lemma\label{root-lemma} For any $n\ge0$ and any real non-singular
homogeneous cubic polynomial $f(x)=f(x_0,\dots,x_n)$, the real
cubic $n$-fold $X\subset P^{n+1}$ defined by equation $z^3=f(x)$
belongs to the deformation component $\Cr$.
\endlemma

\demo{Proof} First, note that $X(\R)$ is diffeomorphic to
$\Rp{n}$, a diffeomorphism being given by projection $P^{n+1}\to
P^n$ forgetting the $z$-coordinate. Moreover, the same projection
provides a homeomorphism between $X(\R)$ and $\Rp{n}$ even if we
authorize $f=0$ to have arbitrary singularities. Combining these
observations with Lemma \ref{through-cusp}, we conclude that the
cubic $z^3=f_t(x)$ does not change the deformation component, when
$f_t$ crosses a facet of the discriminant. To see that this unique
deformation component is $\Cr$, it is sufficient to apply Lemma
\ref{through-cusp} to, say,
$f=x_0(x_1^2+\dots+x_n^2)+x_1^3+\dots+x_n^3$; then one of the two
perturbations described in Lemma \ref{through-cusp} gives
$\Rp{n}\+S^n$, while the other contains elements from $\Cr$.
\qed\enddemo

\proposition\label{many-facets} The boundary of the component
$\Cr$ contains certain facets $\Cal F_q$ with $0\le q\le n$ such
that: $q$ is the index of the facet $\Cal F_q$; the facet $\Cal
F_q$ is adjacent to $\Cal F_p$ through a cuspidal stratum if
$p+q=n$; and the component $C_q=C_p$ adjacent to $\Cr$ through
$\Cal F_q$ and $\Cal F_p$ is formed by cubics $X$ with $X(\R)$
diffeomorphic to $\Rp{n}\#(S^p\times S^q)$.
\endproposition

\demo{Proof} If $f(x)=0$ is a nodal $(n-1)$-fold cubic, the
$n$-fold cubic $z^3=f(x)$ belongs to a cuspidal stratum of
$\D(\R)\subset P_{n,3}(\R)$. It remains, therefore, to apply Lemma
\ref{through-cusp} selecting a suitable polynomial $f$ for each
pair $p, q$. \qed\enddemo

The same argument can be used to construct non-singular real cubic
$n$-folds whose real part is diffeomorphic to $\Rp{n}$ with
several handles.

\proposition For any integer $k$, $0\le k <\frac12(n+1)$, there
exists a non-singular real cubic $n$-fold $X\subset P^{n+1}$ with
$X(\R)$ diffeomorphic to $\Rp{n}\#\mm_k(S^k\times S^{n-k})$, where
$\mm_k$ is equal to  $\binom{n+1}k$. Such a cubic $X$ can be
obtained by perturbation of a certain real cubic $X_0$ with
$\mm_k$ real cusps of index $(k,n-k)$ and with $X_0(\R)=\Rp{n}$.
\endproposition

\demo{Proof} Let us start with a real cubic $(n-1)$-fold defined
in $\Rp{n}$ by a cubic polynomial
$f_t(x)=t(x_0^3+\dots+x_n^3)-(x_0+\dots+x_n)^3$. For $t=(n+1-
2k)^2$ with $0\le k <\frac12(n+1)$, such a cubic has $\mm_k$ nodes
of index $(k,n-k)$.

The real cubic $n$-fold $X_0$ defined in $\Rp{n+1}$ by the
suspension polynomial $f_t(x)-z^3$ has cusps at the corresponding
points, and $X_0(\R)=\Rp{n}$. As in Lemma \ref{through-cusp}, the
perturbation $f_t(x)-z^3-bz-c$ provides, for suitable $b$ and $c$,
a real cubic $n$-fold $X$ with $X(\R)=\Rp{n}\# \mm_k(S^k\times
S^{n-k})$. \qed\enddemo

A different kind of examples is generated by cubics with two nodes
of different indices. Namely, consider $P^{n+2}$ with homogeneous
coordinates $x=(x_0,\dots,x_n)$, $y_1,y_2$, and a cubic $V$
defined by $F(x,y)=f(x)+y_1g_1(x)+y_2g_2(x)+y_1y_2h(x)$, where
$g_i$ are quadratic forms and $f$, $h$ are cubic and linear forms
respectively. In the chart $y_1\ne0$ centered at $(x,y_2)=(0,0)$,
 the affine equation of $V$ is $y_2h+g_1+(f+y_2g_2)$, and, thus, we have a
 node
at the origin if the quadratic form $y_2h+g_1$ is non-degenerate.
 Similarly, in the chart $y_2\ne0$ centered at
$(x,y_1)=(0,0)$ we have a double point governed by the quadratic
form $y_1h+g_2$. It is easy to choose $g_i$ and $h$ so that the
two nodes have arbitrary prescribed indices in the range between
$1$ and $n+1$.
 Then $F(x,y)-z^3$ defines a cubic with two
 cusps which can be perturbed into handles of
 the corresponding indices,
 so that we obtain the following result.

\proposition For any $n\ge0$ and $1\le a,b\le n+1$, there exists
 a real non-singular cubic $(n+2)$-fold, $X\subset P^{n+3}$, with
$X(\R)$ diffeomorphic to $\Rp{n+2}\#(S^a\times
S^{n+2-a})\#(S^b\times S^{n+2-b})$. \qed\endproposition

\subsection{Homological and topological types of
higher dimensional cubics} An important consequence of deformation
classification of real cubic fourfolds in \cite{FK1} is their {\it
homological quasi-simplicity}. By definition, this means that the
triple $(\M(X),c_X,h_X)$ is a complete invariant of a non-singular
cubic fourfold $X$ up to coarse deformation equivalence. The same
is true for $n$-dimensional cubics $X$ if $n<4$. However, in
higher dimensions it is no longer true. We show it below
restricting ourselves to the case of even $n$: in this case we do
not need to involve any techniques essentially different from what
is already discussed in this paper.

The only new ingredient we use in the proof below is a property of
algebraic varieties of type I known as {\it Klein's principle}
(see \cite{Klein}): ``varieties of type I do not admit
development''. According to a modern interpretation (see
\cite{Rokhlin} and \cite{ViroSixyears}), this refers to the local
maximality of such varieties, namely: {\fib a real algebraic
variety of type I cannot increase its total Betti number $b_*$
after crossing the wall in the parameter space}. In the case of
cubics of even dimension, the proof is easy and can be based on
two observations: first, $d=\frac12\big(b_*(X)-b_*(X(\R))\big)$
where, as usual, $d$ is the rank of the 2-periodic {\it
discriminant group} $\M/(\M_++\M_-)=\M_\pm^*/\M_\pm$; second, to
decrease $d$ the vanishing cycle $v$ of the wall crossing should
be of the form $x\pm cx$, which is incompatible with $v^2=\pm 2$
and $x\cdot cx + x^2=0\mod 2$ (the latter is the definition of
type I, see \ref{mark}).

 \proposition For any even $n
\ge6$ there exists a pair of real non-singular cubic $n$-folds
$X,X'\subset P^{n+1}$, $i=1,2$, which are homologically equivalent
but have non-homeomorphic real point sets {\rm (}and thus, $X$ is
not coarse deformation equivalent to $X'${\rm )}.
\endproposition

\demo{Proof} Consider representatives $X_p$ of the deformation
components $C_p$ provided by Proposition \ref{many-facets},
$p=0,\dots,n$. The parity of $p$ determine the rank of each of the
lattices $\M_\pm$, and thus, determine their signature (since one
of the inertia indices is fixed). The discriminant forms have the
same rank (it is less by one than the maximal rank achieved for
$\Cr$) and the same Arf-Brown invariant. Calculation shows that
these lattices are indefinite, and thus, according to Nikulin's
results (see \cite{Nikulin}), there may be only two isomorphism
classes, which correspond to cubics of type I and II,
respectively.

On the other hand, our construction of
$X(\R)=\Rp{n}\#\mm_1(S^1\times S^{n-1})$ and
$X(\R)=\Rp{n}\#\mm_3(S^3\times S^{n-3})$, where $\mm_1=n+1$ and
$\mm_3=\binom{n+1}3$, shows that the cubics $X_p$ admit a
``development'', so according to Klein's principle they are all of
type II.

Now we may conclude that, for example, the cubics $X_1$ and $X_3$
have isomorphic eigenlattices, hence they are homologically
equivalent. But at the same time, for $n\ge6$,
 $X_1(\R)=\Rp{n}\#(S^1\times S^{n-1})$ and $X_3=\Rp{n}\#(S^3\times S^{n-3})$
are not homeomorphic.\qed
\enddemo

\subsection{
Higher dimensional cubics as ramified connected sums} According to
Lemma \ref{nonsingularity}, whatever is the dimension, any family
of non-singular cubics $X_t$ degenerating to a nodal cubic $X_0$
can be transformed into a family of cubics $X_t'$, such that $X_t$
is projectively equivalent to $X_t$, and $X_t'$ degenerates to a
cubic splitting into a union of a quadric and a hyperplane. This
remarkable property of cubics implies that all the deformation
classes of non-singular real cubics can be obtained by small
perturbation of such reducible cubics. More precisely, the
following statement holds (and follows easily from Lemma
\ref{nonsingularity} and the possibility to degenerate any
non-singular cubic to a nodal one).

\proposition Any deformation class of non-singular cubic
hypersurfaces in $P^{n+1}$, with coordinates $x_0,\dots,x_{n},y$,
admits a representative defined by equation $\{y(f_2-y^2)+\delta
f_3=0\}$, where $f_2=f_2(x)=\sum_{i=0}^{n}\pm x_i^2$ (for a
certain choices of signs $\pm$), $0<\delta<\!<1$, and
$f_3=f_3(x,y)$ is some cubic transversal to the quadric $f_2$ in
the hyperplane $y=0$. \qed\endproposition

Combining this statement with Proposition \ref{rsum}, we obtain
the following description of topology of real loci.

\corollary The real point set $X(\R)$ of any non-singular real
cubic hypersurface $X\subset P^{n+1}$, for any $n>0$, is
diffeomorphic to a ramified connected sum $P\#^FD(P_+)$, where
$P=P^{n}(\R)$, $D(P_+)$ is the double of $P_+=\{x\in
P\,|\,f_2\ge0\}$, and $F=\{x\in P_+\,|\,f_3=0\}$. Here $f_2$,
$f_3$ are respectively quadratic and cubic forms in $P$, such that
$f_2$
 is non-degenerate and has a transverse intersection with $f_3$.
\qed\endcorollary


\Refs\widestnumber\key{ABC}

\ref{ACT}
 \by D. Allcock, J. Carlson, and D. Toledo
 \paper Hyperbolic geometry and moduli of real cubic surfaces
 \jour arXiv:0707.1058 \toappear \ in Ann. Sci. ENS
 \vol
 \issue
 \yr
 \pages
\endref\label{ACT}

\ref{DK} \paper  Topological properties of real algebraic
varieties : du c\^{o}t\`{e} de chez Rokhlin \by A.~Degtyarev and
V.~Kharlamov \jour Uspekhi Mat. Nauk. \vol 55 \yr 2000 \issue 4
\pages  129--212
\endref\label{DK}

\ref{FK1}
 \by S. Finashin, V. Kharlamov
 \paper Deformation classes of real four-dimensional cubic
 hypersurfaces
 \jour J. Alg. Geom.
 \vol17
 \issue
 \yr2008
 \pages677--707
\endref\label{FK1}

\ref{FK2}
 \by S. Finashin, V. Kharlamov
 \paper On the deformation chirality of real cubic
fourfolds
 \jour arXiv:0804.4882
\endref\label{FK2}

\ref{Ki} \by R.~Kirby \paper A Calculus for Framed Links in $S^3$
\jour Invent. Math. \vol 45 \yr 1978 \pages 35--56
\endref\label{Kirby}


\ref{Kl}
 \by F.~Klein
 \book Gesammelte mathematische Abhandlungen
 \bookinfo Berlin
 \vol 2
 \yr 1922
 \endref\label{Klein}

\ref{Kr1}
 \by V.~Krasnov
 \paper Rigid isotopy classification of real three-dimensional
 cubics
 \jour Izvestiya: Mathematics
 \vol 70
 \issue 4
 \yr2006
 \pages 731--768
\endref\label{Kr1}

\ref{Kr2}
 \by V.~Krasnov
\paper On the topological classification of three-dimensional real
cubics \jour \toappear \ in Izvestiya: Mathematics
\endref\label{Kr2}



\ref{La}
 \by R. Laza
 \paper The moduli space of cubic fourfolds via the period map
 \jour arXiv:0705.0949,
\toappear \ in Annals of Math
 \vol
 \issue
 \yr
 \pages
\endref\label{Laza}

\ref{Lo}
 \by E. Looijenga
 \paper The period map for cubic fourfolds
 \jour arXiv:0705.0951, \toappear \ in Invent. Math.
 \vol
 \issue
 \yr
 \pages
\endref\label{Looijenga}

\ref{N}
 \by V.~V.~Nikulin
 \paper Integer quadratic forms and some of their geometrical applications
 \jour Math. USSR -- Izv.
 \vol 43
  \yr 1979
 \pages 103--167
\endref\label{Nikulin}

\ref{O}
 \by P. Orlik
\book Seifert manifolds, \bookinfo Lecture Notes in Math.,
Springer-Ver., Berlin - New York
 \vol291
 \yr1972
\endref\label{Orlik}

\ref{R}
 \by V.~A.~Rohlin
 \paper Complex topological characteristics of real algebraic curves
 \jour Uspekhi Mat. Nauk
 \vol 33
  \yr 1978
 \pages 77--89
\endref\label{Rokhlin}


\ref{Vi}
 \by O.~Viro
 \paper Progress of the last six years in topology of real algebraic varieties
 \jour Russian Math. Surveys
 \vol 41
  \yr 1986
 \pages 55--82
\endref\label{ViroSixyears}

\ref{Vo} \by C.~Voisin
 \paper Th\'{e}or\`{e}me de Torelli pour les cubiques de $P^5$
 \jour Invent. Math.
 \vol 86
 \yr 1986
 \pages 577--601
\endref\label{Voisin}


\endRefs
\enddocument